\documentclass[a4,11pt]{article}
\usepackage{amsmath,amsthm,amssymb}
\usepackage{graphicx}

\title{}
\author{}
\date{}

\newtheorem{Th}{Theorem}[subsection] 
\newtheorem{Lem}[Th]{Lemma} 
\newtheorem{Prop}[Th]{Proposition} 
\newtheorem{Cor}[Th]{Corollary} 
\newtheorem{Def}[Th]{Definition} 
 
\newtheorem{Assum}[Th]{Assumption} 
\newtheorem{Ex}[Th]{Example}

\newcommand{\ba}{\mbox{\boldmath $a$}} 
\newcommand{\sba}{\mbox{\scriptsize \boldmath $a$}} 
 
\newcommand{\bx}{\mbox{\boldmath $x$}} 
\newcommand{\balpha}{\mbox{\boldmath $\alpha$}} 
\newcommand{\sbx}{\mbox{\scriptsize \boldmath $x$}} 
\newcommand{\by}{\mbox{\boldmath $y$}} 
\newcommand{\sby}{\mbox{\scriptsize \boldmath $y$}}
 
\newcommand{\bnu}{\mbox{\boldmath $\nu$}}

\def\N{\mathbb{N}} 
\def\R{\mathbb{R}}

\def\D{{\cal D}}

\def\K{{\cal K}}

\def\cI{{\cal I}} 

\def\DS{\displaystyle} 
\def\1/2{\frac{1}{2}}

\def\dist{\mbox{\rm dist}} 
 
\def\div{\mbox{\rm div}} 
\def\supp{\mbox{\rm supp}}

\def\ov{\overline} 
 
\def\vep{\varepsilon} 
 
\def\1/2{\frac{1}{2}}

\def\esup{\sup}
\def\einf{\inf}
\def\rT{{\rm T}}
\def\loc{{\rm\scriptsize loc}}

\makeatletter 
\long\def\@makefntext#1{\parindent 1em\noindent 
\@hangfrom{\hbox to 1.8em{\hss$^{\@thefnmark}$}}#1}
\makeatother

\begin{document}
\mbox{}
\vspace{-0.5cm}
\begin{center}
{\large\bf Exponential decay phenomenon of the principal eigenvalue\\
of an elliptic operator with a large drift term of gradient type}\\
~\\
{Shuichi Jimbo$^*$, Masato Kimura$^{**}$, Hirofumi Notsu$^{**}$}\\
~\\
{\small
\begin{tabular}{rl}
$^{\rm *}$ &
{\it Department of mathematics, Hokkaido University,
Sapporo, 060-0810, Japan} \\ \vspace{1ex}
~&{\it E-mail: jimbo@math.sci.hokudai.ac.jp}\\ 
$^{\rm **}$ &
{\it Faculty of Mathematics, Kyushu University,
Fukuoka, 812-8581, Japan}\\
~&{\it E-mails: \{masato,~notsu\}@math.kyushu-u.ac.jp}
\end{tabular}
}
\end{center}
\vspace{2ex}

\begin{abstract} 
We study an asymptotic behaviour of the principal eigenvalue for an elliptic 
operator with large advection which is given by a gradient of
a potential function. 
It is shown 
that the principal eigenvalue decays exponentially 
under the velocity potential well condition as the 
parameter tends to infinity.
We reveal that the depth of the potential well
plays an important role in the estimate.
Particularly, in one dimensional case,
we give a much more elaborate characterization for the
eigenvalue.
Some numerical examples obtained by a characteristic-curve finite element 
method are also shown. 
\end{abstract} 

\subsection{\large Introduction}\label{intro} 
\setcounter{equation}{0}

The following elliptic eigenvalue problem with a large drift term 
is considered in this paper.
Let $\Omega$ be a bounded Lipschitz domain in $\R^n$ ($n\geq 1$).
For a given vector field 
$\ba =(a_1,\cdots,a_n)^\rT 
\in L^\infty (\Omega,\R^n)$, we consider 
an elliptic eigenvalue problem with a parameter $p\in\R$ :
\begin{equation}\label{evp}
\left\{
\begin{array}{ll}
-\Delta u (\bx)+p\,\ba(\bx)\cdot\nabla u (\bx)=\lambda u (\bx)
&(\bx\in\Omega)\\~\\
u (\bx)=0
&(\bx\in\partial \Omega)
\end{array}
\right.
\end{equation}
We assume that the eigenfunction $u$ ($u\not\equiv 0$)
belongs to the Sobolev
space $H^1_0(\Omega)$. Then, from the elliptic regularity
theorem (see \cite{G-T83} for example), 
$u\in W^{2,s}_\loc (\Omega)$ holds
for arbitrary $s\in [1,\infty)$.

Although $\lambda$ and $u (\bx)$ are complex-valued
in general, we can define a positive real-valued principal 
eigenvalue $\lambda=\lambda_1(p)$ and the corresponding
real-valued eigenfunction $u=u_1(\cdot,p)$,
which is uniquely determined by the condition:
\begin{equation}\label{ucondition}
u (\bx)>0~~~(\bx\in\Omega),
~~~
\max_{\sbx\in\Omega}u (\bx)=1\,.
\end{equation}
Some related results for the principal eigenvalue
and eigenfunction will be collected in 
Theorem~\ref{fth} and Theorem~\ref{ONS2}.
A detailed and systematic study for the
principal eigenvalue for general second order 
elliptic operators in general domains
can be found in 
Berestycki, Nirenberg and Varadhan \cite{B-N-V94}.

As we will see in Section \ref{numerical} through several
numerical examples, the principal eigenvalue
is closely related to the decay rate of 
the solution of the corresponding nonstationary linear advection-diffusion
equation (\ref{parabolic}).

Our interest in this study is
asymptotic behaviours of $\lambda_1(p)$
as $p\rightarrow \infty$. 
Such large effects from advection or drift term
under small diffusion
appear in many actual physical problems
and often cause some difficulties in 
their numerical simulations and analysis.

This phenomenon was first studied by 
Ventcel' \cite{Ven72} and Friedman \cite{Fri73} 
with Ventsel'-Freidlin's probabilistic 
approach \cite{V-F70}. 
In \cite{Fri73}, the following result
is obtained.
If $\Omega$ is a bounded domain of $C^2$-class,
and if $\ba$ belongs to $C^1(\ov{\Omega})$ and
satisfies the condition:
\begin{equation}\label{Fcondition}
\ba(\bx)\cdot\bnu(\bx)>0
~~~~~(\bx\in\partial\Omega),
\end{equation}
where $\bnu$ denotes the outward unit normal vector
on $\partial \Omega$,
there exist $c_1$ and $c_2$ such that
\[
0<c_1\leq \liminf_{p\to\infty}
\frac{1}{p}\log \frac{1}{\lambda_1(p)}
\leq
\limsup_{p\to\infty}
\frac{1}{p}\log \frac{1}{\lambda_1(p)}
\leq c_2.
\]
This estimate is equivalent to 
\[
^\forall \vep >0,
~~~
^\exists p_0\in\R
~~~
\mbox{s.t.}
~~~
e^{-(c_2+\vep)p}\leq \lambda_1(p)
\leq e^{-(c_1-\vep)p}
~~~~~
(^\forall p\geq p_0),
\]
which means that the principal eigenvalue
becomes exponentially small under the condition
(\ref{Fcondition}).

In one dimensional case, asymptotic behaviours of $k$th 
eigenvalues $\lambda_k(p)$ $(k\in\N)$
are studied in \cite{deG80a} by means of matched 
asymptotic expansions. In connection with a singular
limit analysis of Sturm-Liouville two points boundary 
value problems, precise approximations
of eigenvalues and eigenfunctions including exponentially small
principal eigenvalue are obtained there.

For a given vector field $\ba\in L^\infty (\Omega,\R^n)$,
if the principal eigenvalue $\lambda_1(p)$ satisfies
the condition:
\begin{equation}\label{edp1}
\mbox{
$^\exists c>0$ and $^\exists p_0\in\R$ ~~s.t.~~~}
0<\lambda_1(p)\leq e^{-cp}~~~
(^\forall p\geq p_0),
\end{equation}
or its equivalent condition:
\[
\liminf_{p\to\infty}
\frac{1}{p}\log \frac{1}{\lambda_1(p)}>0,
\]
then we call it 
exponential decay phenomenon of 
the principal eigenvalue.
In addition to some numerical
examples in Section~\ref{numerical},
we will give a biological interpretation 
in a chemotaxis model
for more physical pictures of the
exponential decay phenomenon.

Besides the exponential decay phenomenon,
the principal eigenvalue $\lambda_1(p)$ 
exhibits various asymptotic behaviour 
as $p\to\infty$.
Devinatz, Ellis and Friedman \cite{D-E-F74} investigated 
$L^2$ inequalities type arguments
and maximum principle type arguments,
and they obtained some estimates of
$\lambda_1(p)$ from above or from below
by a term $Cp^\gamma$ with $\gamma\in (0,2]$.
See also \cite{Ven75} for related results.

In \cite{B-H-N05}, Berestycki, Hamel and Nadirashvili
proved that, if $\partial\Omega$ is of $C^2$-class 
and $\ba\in L^\infty (\Omega,\R^n)$ satisfies
$\div \ba =0$ in $\D'(\Omega)$, then
$\limsup_{p\to\infty} \lambda_1(p)<\infty$ if and only if
there exists $\psi\in H^1(\Omega)$ such that
$\psi\not\equiv \mbox{const}$ and $\ba\cdot\nabla\psi=0$
a.e.\,in $\Omega$.

In our study, we focus on the exponential decay phenomenon
and give some new estimates for it. 
We assume existence of a velocity potential of $\ba$
(condition (\ref{assumeb})),
which allows us to use the Rayleigh quotient 
and $L^2$ inequalities type arguments.
One of our main results shows us that
existence of a potential well (Definition~\ref{pwell})
implies the exponential decay phenomenon of 
the principal eigenvalue, and that
the depth of the potential well
gives an estimate of the constant $c$ in (\ref{edp1})
from below (Theorem~\ref{mth}).
This observation is extended to the exponential decay phenomenon of
the $m$th eigenvalue (Theorem~\ref{mth2})
and to a precise asymptotic behaviour of $\lambda_1(p)$
in one dimensional case (Theorem~\ref{1dth}).  

The organization of this paper is as follows.
In Section~\ref{numerical}, 
we show some numerical examples of the exponential decay phenomenon
in two dimensional case.
These  numerical profiles of 
the principal eigenfunctions $u_1(\bx,p)$ 
will be helpful in our analysis later.
We collect some fundamental facts on the principal eigenvalue
in Section~\ref{preliminaries}.
In Section~\ref{asymptotic},
we prove several estimates for asymptotic behaviours of
$\lambda_1(p)$ including an alternative simpler proof of the
exponential decay phenomenon.
In Section~\ref{subsec2}, we give more precise asymptotic
behaviour of $\lambda_1(p)$ in one dimensional case,
under a different assumption than one of \cite{deG80a}.
Additionally as an application,
we give a biological interpretation of the
exponential decay phenomena, life span of a
biological colony under the chemotaxis effect
in the last section.

\subsection{\large Numerical examples of exponential decay phenomena}
\label{numerical} 
\setcounter{equation}{0}

For the principal eigenvalues $\lambda_1(p)$ 
which is defined in (\ref{evp}),
we give several examples of 
typical asymptotic behaviours of $\lambda_1(p)$
as $p\to\infty$ for fixed vector fields $\ba(\bx)$,
particularly focusing our interests on the exponential decay phenomenon.
We start from the simplest case that $\ba$ is a constant 
vector field. 
\begin{Ex}{\rm
Let $\ba\in\R^n$ be a constant vector. Then
we can easily check that
\[
\lambda_1 (p)=\lambda_1(0)+p^2 \frac{|\ba|^2}{4},
~~~~~
u_1(\bx,p)=e^{\frac{p}{2}\sba\cdot\sbx}\,u_1(\bx,0),
\]
where $\lambda_1(0)$ and $u_1(\bx,0)$ are the 
principal eigenvalue of $-\Delta$ with the Dirichlet boundary
condition on $\Omega$ and its corresponding 
eigenfunction.
Moreover, not only the principal eigenvalue
but also all eigenvalues are exactly shifted by $(|\ba|^2/4)\,p^2$,
i.e., $\lambda_k (p)=\lambda_k(0)+(|\ba|^2/4)\,p^2$ for all $k\in\N$.
In this case, the asymptotic behaviour of $\lambda_1(p)$ (or $\lambda_k(p)$)
is of $O(p^2)$ as $p\to\infty$.
We will see in the next section 
that this type of $O(p^2)$ behaviour is most common,
for example if the vector field $\ba(\bx)$ has no
singular point. Of course, there is no chance to have
the exponential decay phenomenon of principal eigenvalues
as $p\to\infty$.
}\end{Ex}

On the other hand, it seems to be difficult
to give an example of the exponential decay phenomenon
similarly only by using elementary function.

In this Section,
we give some examples of the exponential decay phenomenon
with the help of numerical simulations.
For this purpose, we compute the following initial-boundary value
problem of parabolic type:
\begin{equation}\label{parabolic}
\left\{
\begin{array}{ll}
\DS{u_t-\Delta u +p\,\ba(\bx)\cdot \nabla u=0}
&(\bx\in \Omega, \ 0<t<\infty)\\~\\
u (\bx,t)=0
&(\bx\in\partial\Omega,~t>0) \\~\\
u (\bx,0)=u_0(\bx)
&(\bx\in\Omega),
\end{array}\right.
\end{equation}
where  $u_0$ is a given
positive initial function.
In the following 
numerical examples, 
the domain is chosen as 
$\Omega = (-1,1)\times (-1,1)$.
It is well-known that 
$u (\bx,t)$ uniformly converges to zero as $t\to\infty$
and its asymptotic behaviour is more precisely given as
\begin{equation}\label{asymptotics}
u (\bx,t) \sim  C e^{- \lambda _1(p)t} u_1(\bx,p)
~~~\mbox{as}~t \to \infty,
\end{equation}
where the constant $C$ is positive if $u_0(\bx)>0$ for $\bx\in\Omega$.

If $p>>1$, since (\ref{parabolic}) 
is related to exponentially small positive eigenvalue
and is a convection-dominant problem,
it requires us to adopt some upwinding technique
for its reliable and reasonable numerical simulation.
The following simulations were computed by means of
a characteristic-curve finite element scheme 
(103) in Chap.3 of \cite{Pir89}
with piecewise linear triangular elements.
This is a upwind implicit scheme of first order 
based on the characteristic curve approximation:
\[
u_t(\bx,k\tau)+p\ba(\bx)\cdot\nabla u (\bx,k\tau)
=\frac{u (\bx,k\tau)-u (\bx-p\ba(\bx)\tau,(k-1)\tau)}{\tau}
+O(\tau),
\]
where $\tau>0$ is a small time increment.
For more details and the stability analysis, 
see \cite{Pir82} and \cite{Pir89} etc. 

We use the following vector field with compact support on $\R^2$;
\[
\balpha (\bx; R) := 
\left\{ 
\begin{array}{cl}
\dfrac{\bx}{|\bx | } \sin\left( \dfrac{\pi |\bx |}{R} \right) &
( 0<|\bx | \le R )\\~\\
\textrm{\boldmath $0$} & (\bx=\textrm{\boldmath $0$},~\mbox{or}~|\bx |> R )
\end{array}\right.
\]
We remark that $\balpha (\, \cdot \, ; R)\in W^{1,\infty}(\R^2;\R^2)$
and $\| \balpha (\cdot,R)\|_{L^\infty(\R^2)} = 1$.

\begin{figure}[htbp]
\begin{center}
\begin{minipage}{0.45\hsize}
\begin{center}
\rotatebox{270}{\includegraphics[width=0.75\hsize,clip]{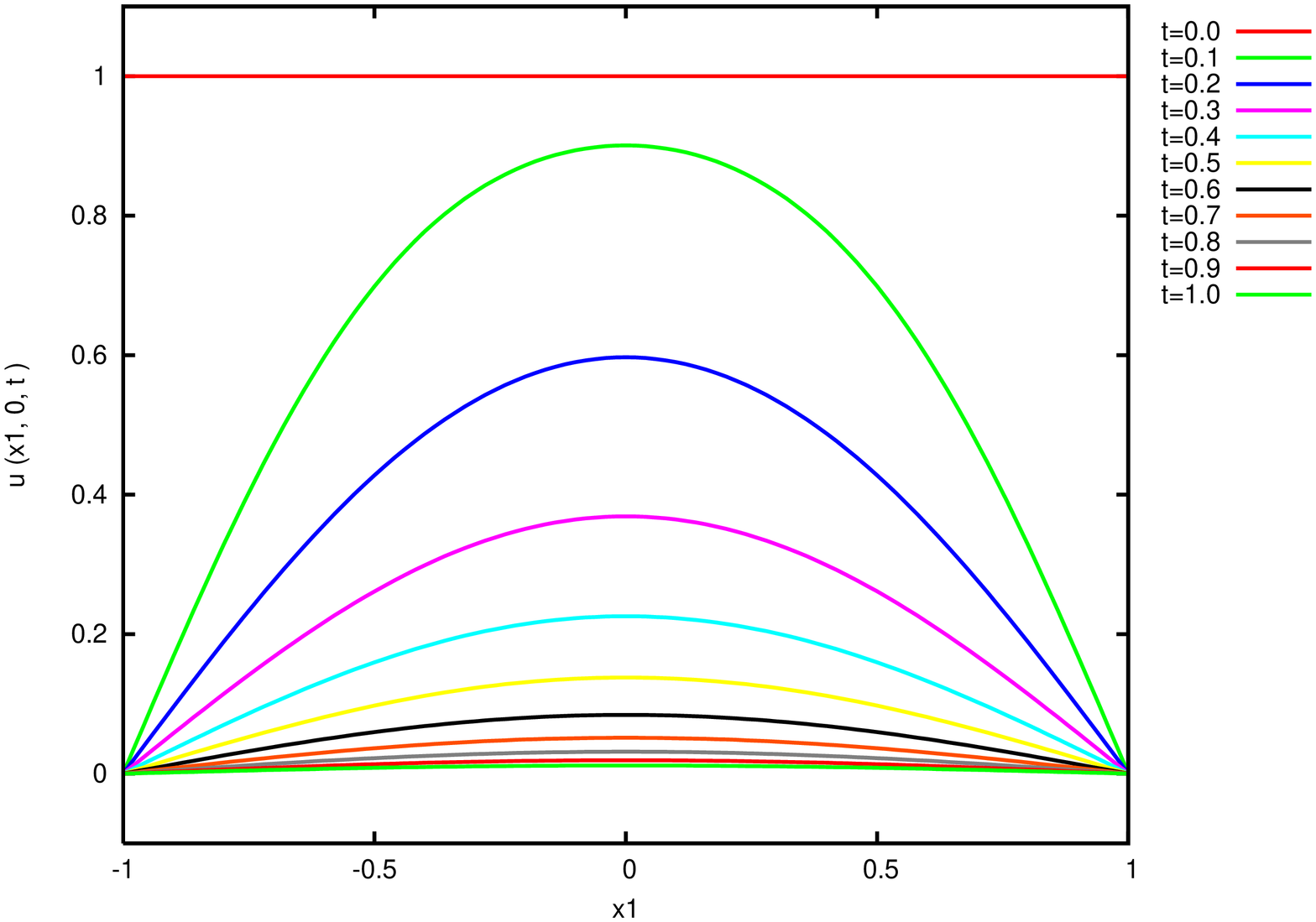}}
\caption{Time evolution of Example~\ref{Sim1} on the line $x_2=0$ for $p=0$.}
\label{Sim1-1d-p00}
\end{center}
\end{minipage}
\hspace{0.08\hsize}
\begin{minipage}{0.45\hsize}
\begin{center}
\rotatebox{270}{\includegraphics[width=0.75\hsize,clip]{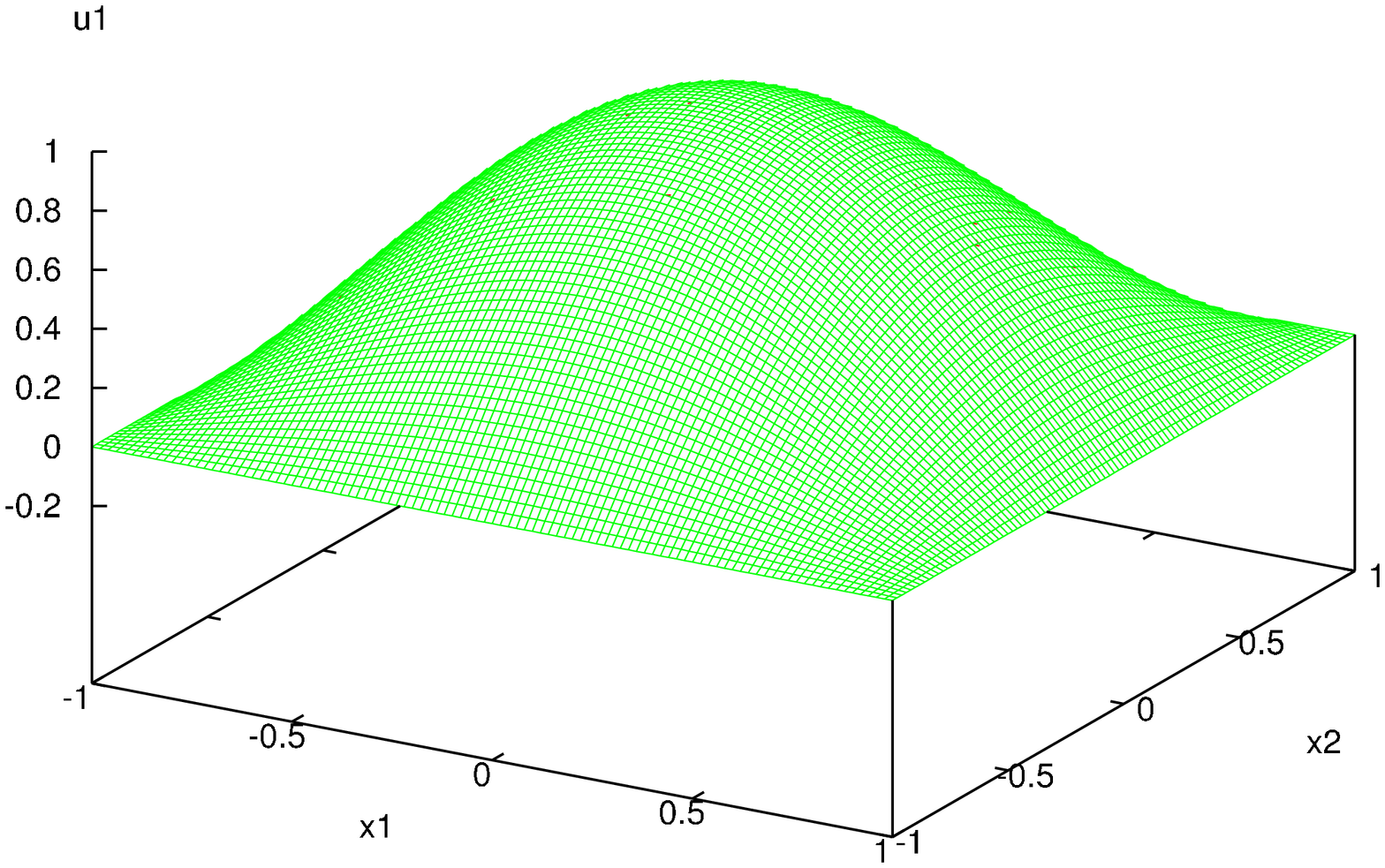}}
\caption{Eigenfunction of Example~\ref{Sim1} for $p=0$.}
\label{Sim1-profile-p00}
\end{center}
\end{minipage}
\end{center}
\begin{center}
\begin{minipage}{0.45\hsize}
\begin{center}
\rotatebox{270}{\includegraphics[width=0.75\hsize,clip]{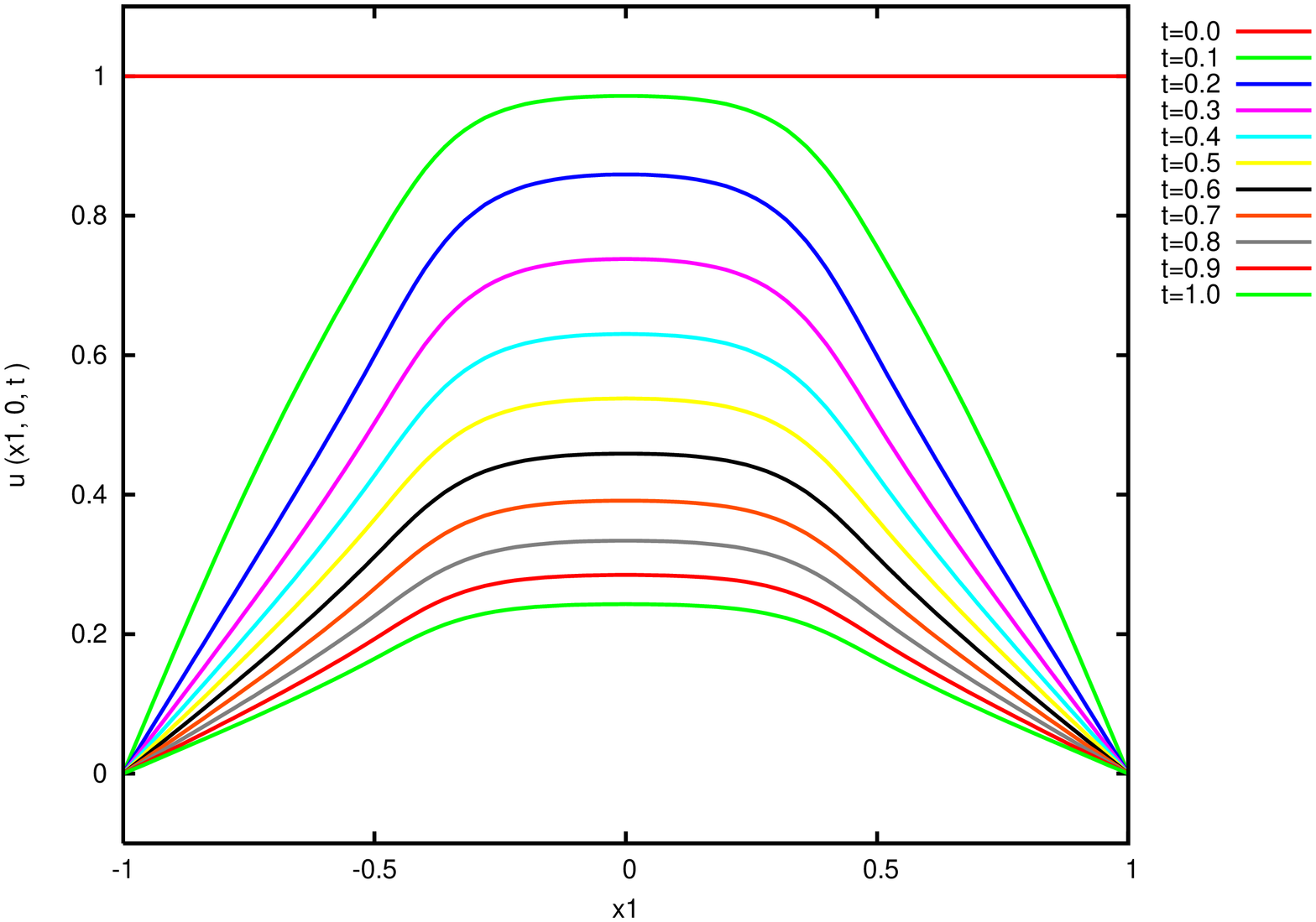}}
\caption{Time evolution of Example~\ref{Sim1} on the line $x_2=0$ for $p=10$.}
\label{Sim1-1d-p10}
\end{center}
\end{minipage}
\hspace{0.08\hsize}
\begin{minipage}{0.45\hsize}
\begin{center}
\rotatebox{270}{\includegraphics[width=0.75\hsize,clip]{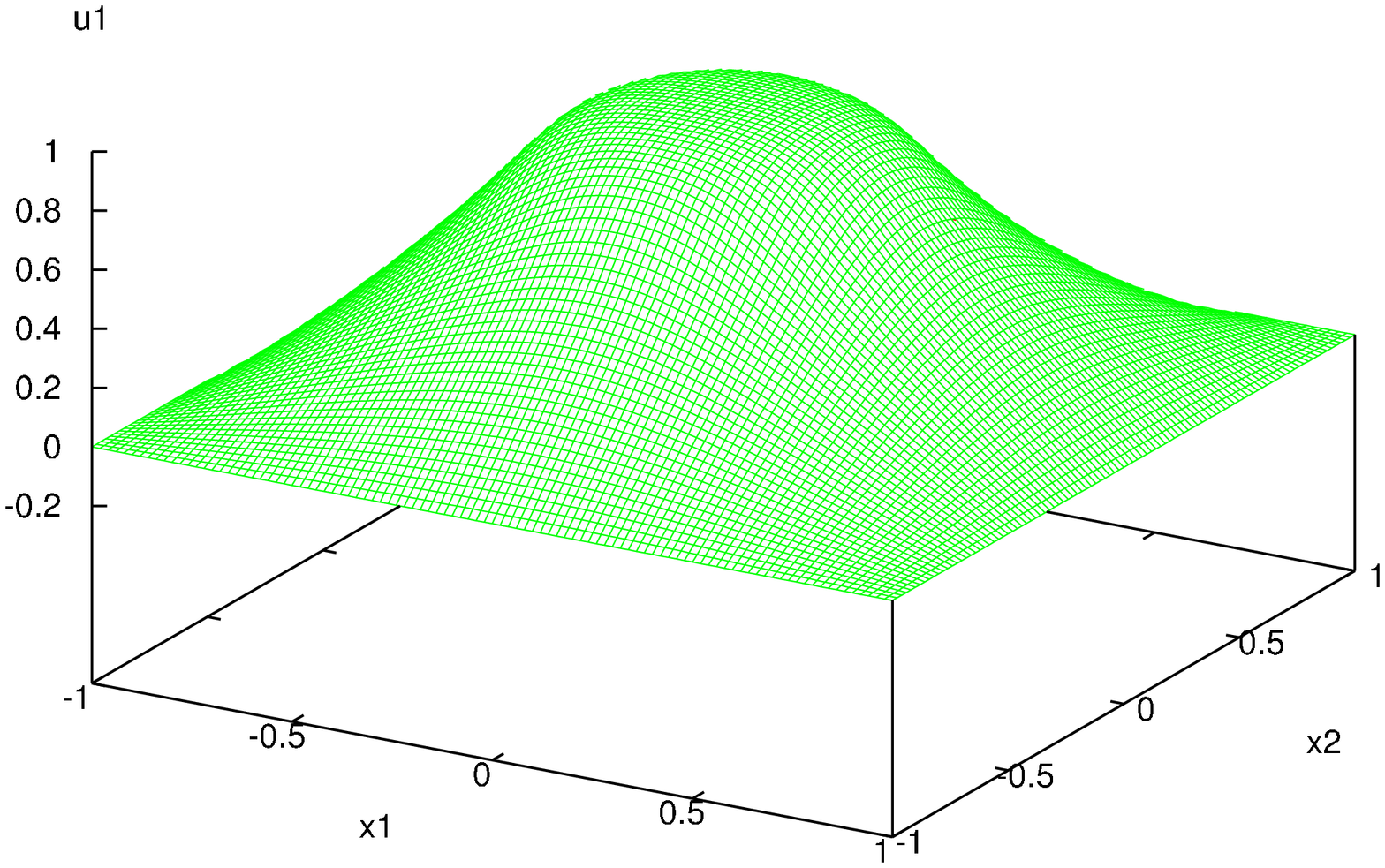}}
\caption{Eigenfunction of Example~\ref{Sim1} for $p=10$.}
\label{Sim1-profile-p10}
\end{center}
\end{minipage}
\end{center}
\begin{center}
\begin{minipage}{0.45\hsize}
\begin{center}
\rotatebox{270}{\includegraphics[width=0.75\hsize,clip]{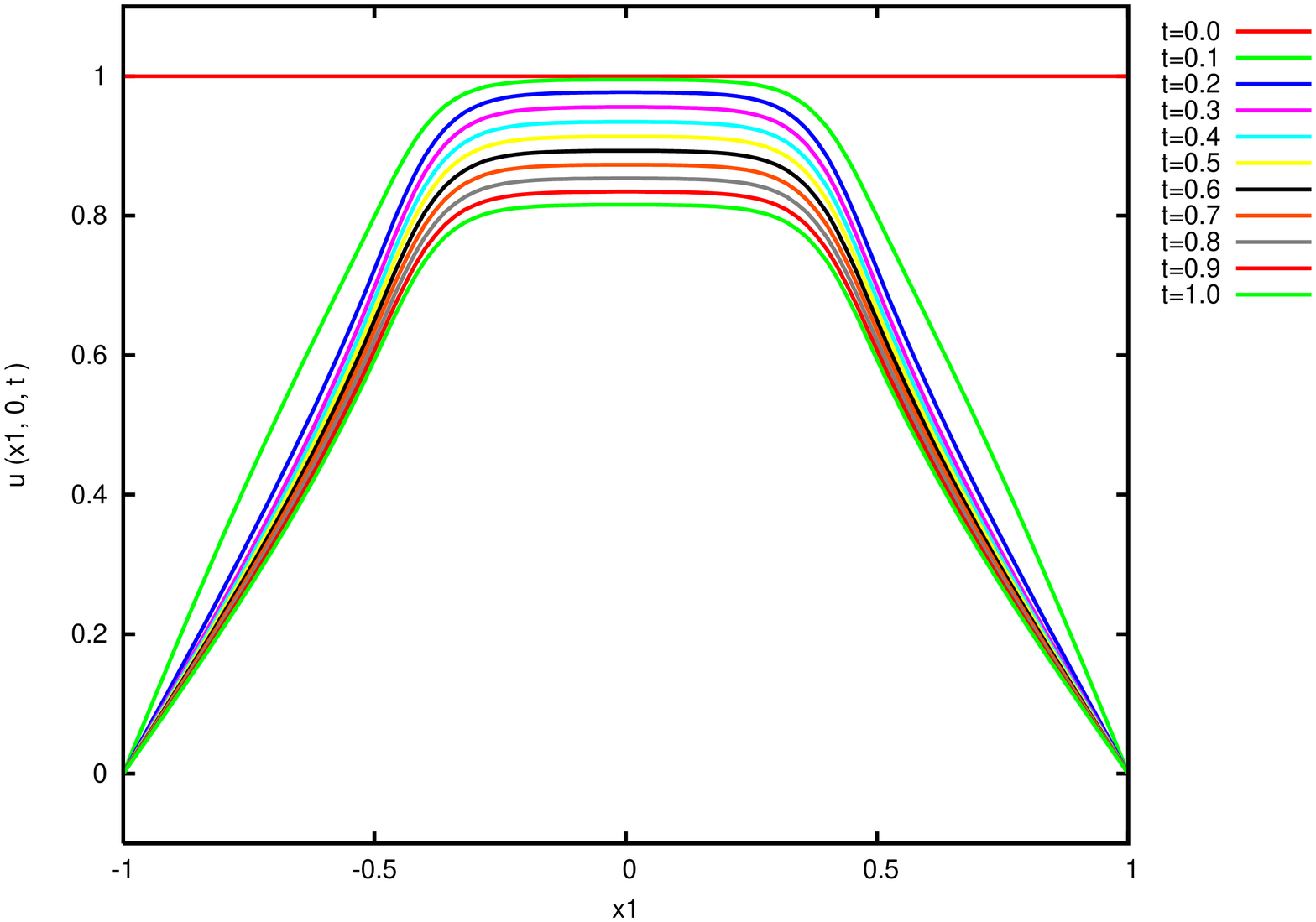}}
\caption{Time evolution of Example~\ref{Sim1} on the line $x_2=0$ for $p=20$.}
\label{Sim1-1d-p20}
\end{center}
\end{minipage}
\hspace{0.08\hsize}
\begin{minipage}{0.45\hsize}
\begin{center}
\rotatebox{270}{\includegraphics[width=0.75\hsize,clip]{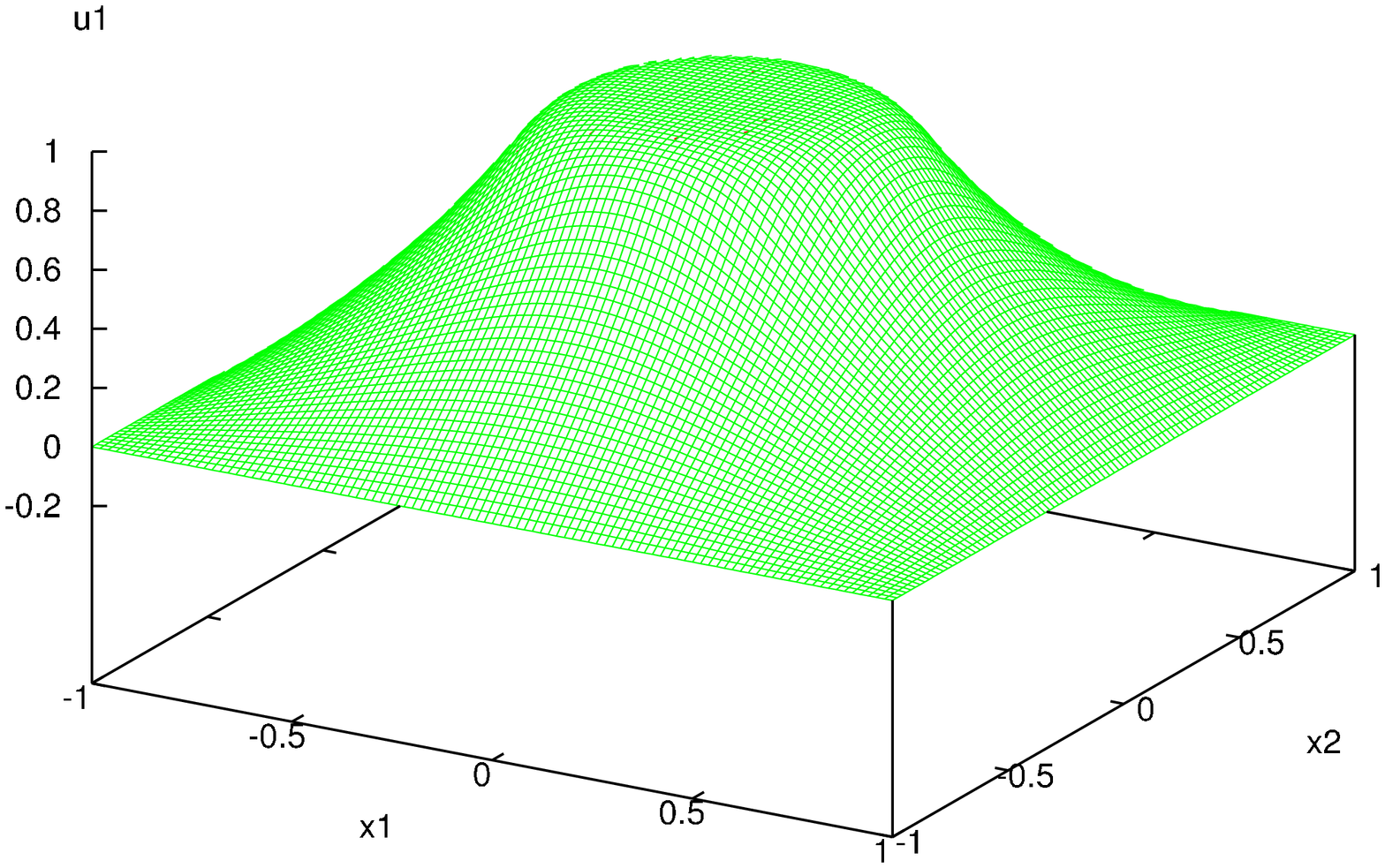}}
\caption{Eigenfunction of Example~\ref{Sim1} for $p=20$.}
\label{Sim1-profile-p20}
\end{center}
\end{minipage}
\end{center}
\end{figure}

\begin{figure}[htbp]
\begin{center}
\begin{minipage}{0.45\hsize}
\begin{center}
\rotatebox{270}{\includegraphics[width=0.75\hsize,clip]{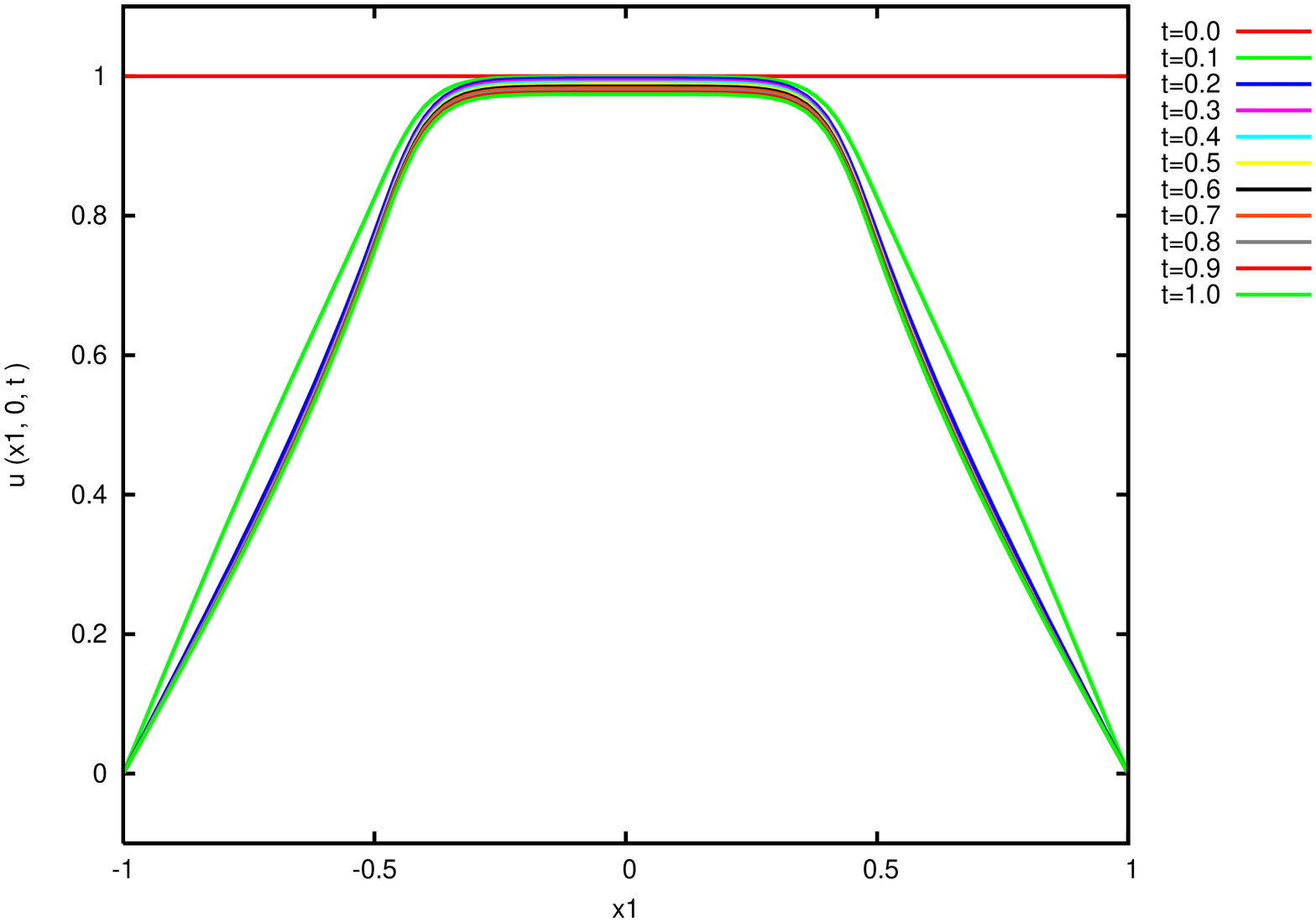}}
\caption{Time evolution of Example~\ref{Sim1} on the line $x_2=0$ for $p=30$.}
\label{Sim1-1d-p30}
\end{center}
\end{minipage}
\hspace{0.08\hsize}
\begin{minipage}{0.45\hsize}
\begin{center}
\rotatebox{270}{\includegraphics[width=0.75\hsize,clip]{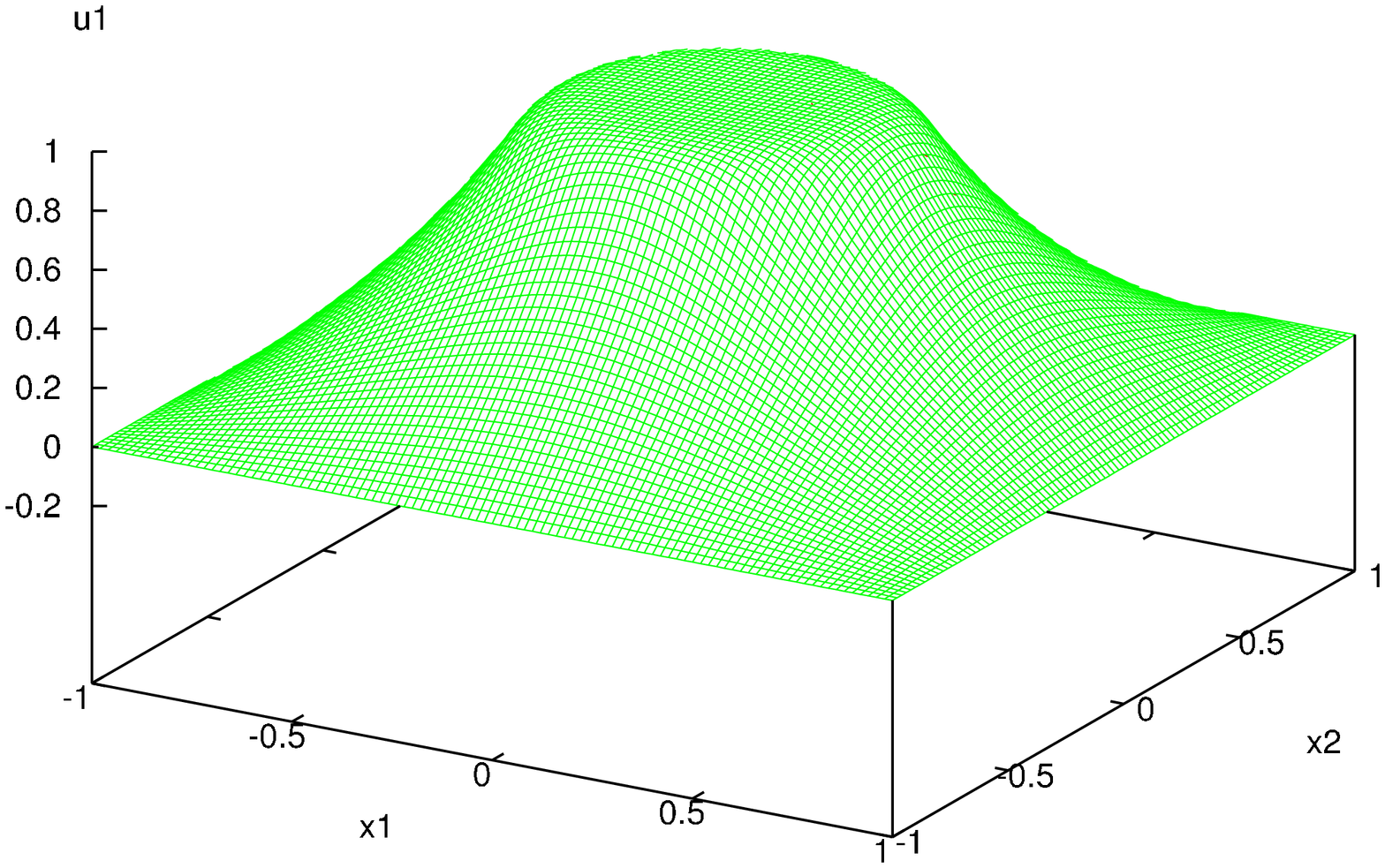}}
\caption{Eigenfunction of Example~\ref{Sim1} for $p=30$.}
\label{Sim1-profile-p30}
\end{center}
\end{minipage}
\end{center}
\begin{center}
\begin{minipage}{0.45\hsize}
\begin{center}
\rotatebox{270}{\includegraphics[width=0.75\hsize,clip]{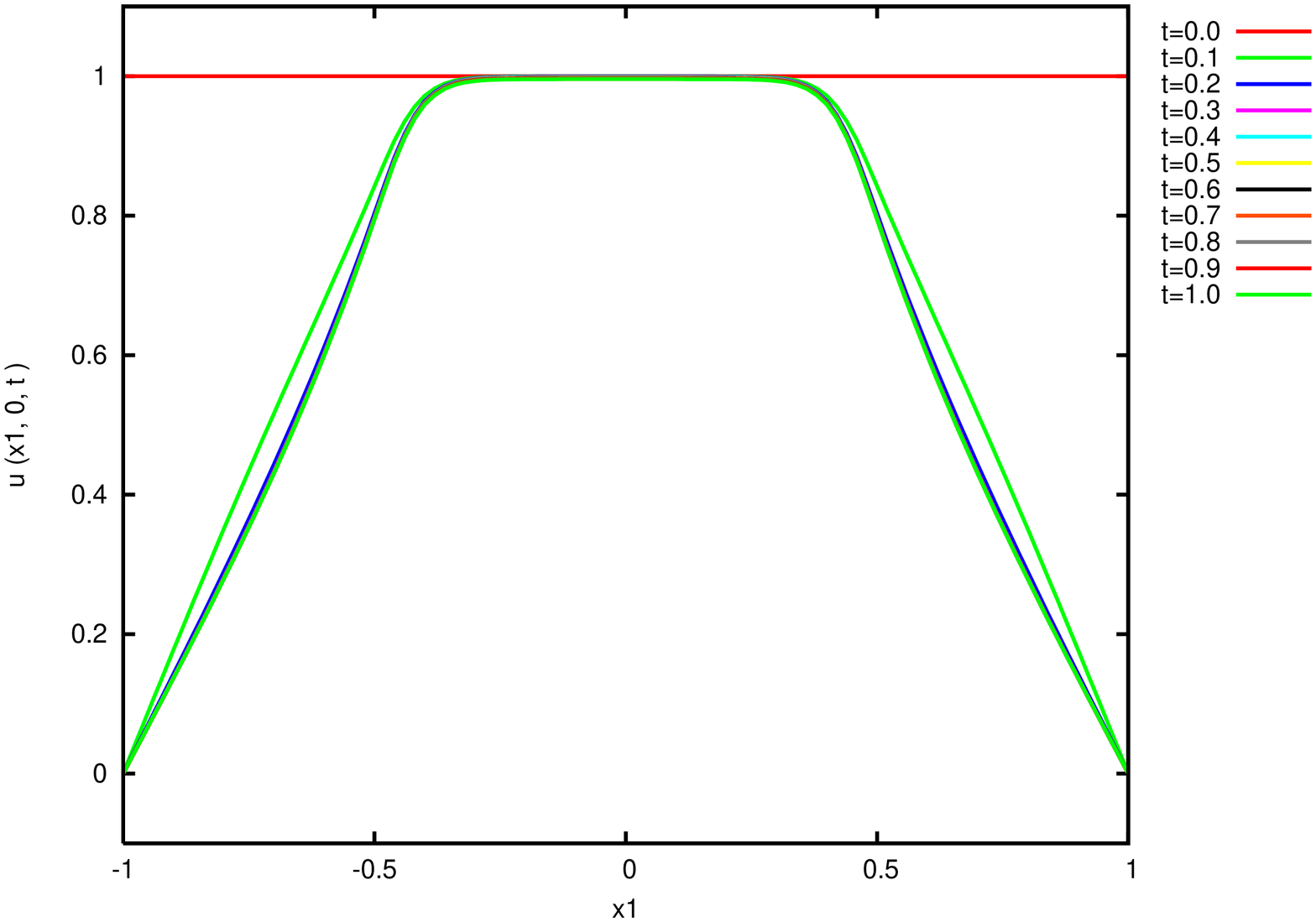}}
\caption{Time evolution of Example~\ref{Sim1} on the line $x_2=0$ for $p=40$.}
\label{Sim1-1d-p40}
\end{center}
\end{minipage}
\hspace{0.08\hsize}
\begin{minipage}{0.45\hsize}
\begin{center}
\rotatebox{270}{\includegraphics[width=0.75\hsize,clip]{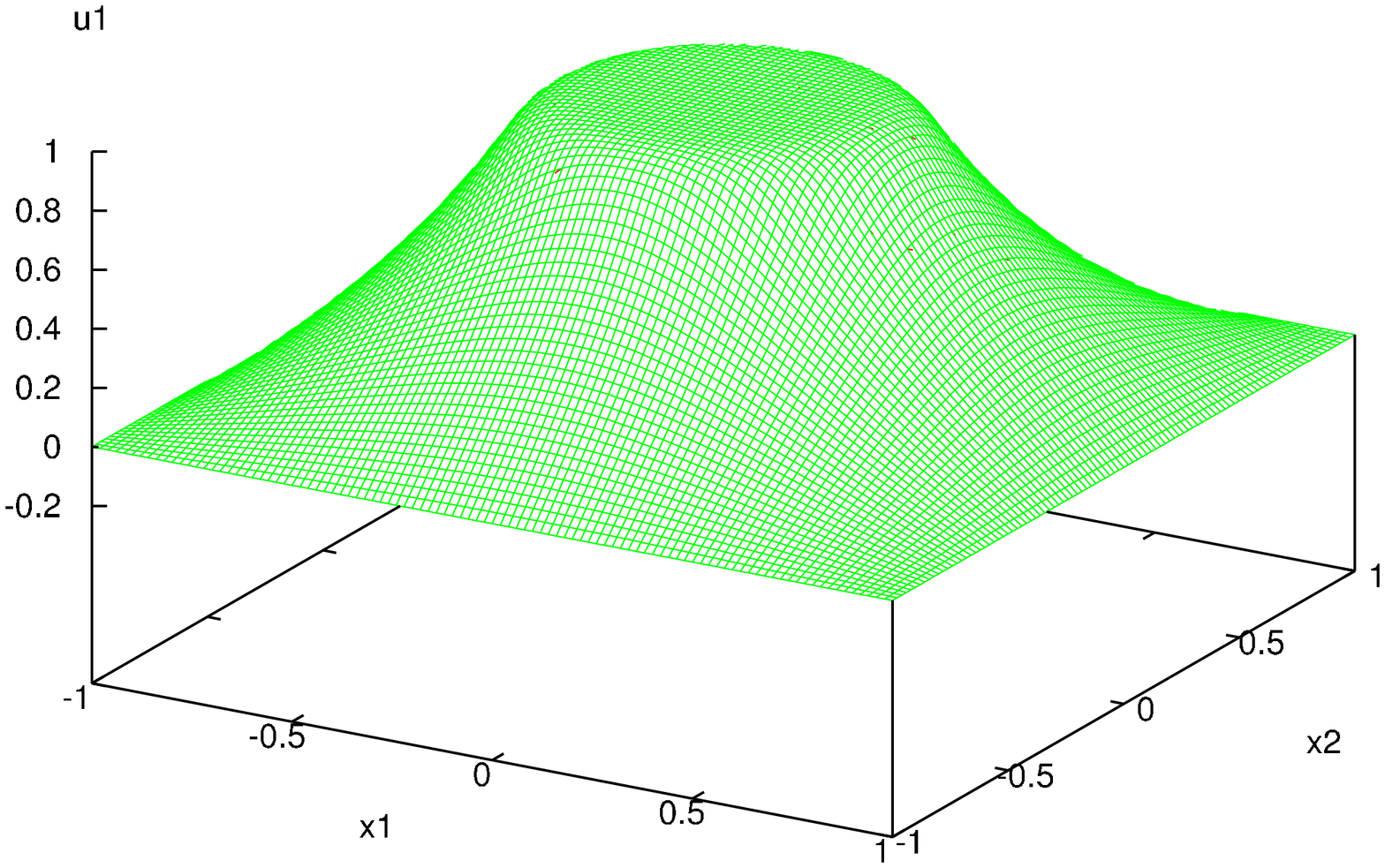}}
\caption{Eigenfunction of Example~\ref{Sim1} for $p=40$.}
\label{Sim1-profile-p40}
\end{center}
\end{minipage}
\end{center}
\begin{center}
\begin{minipage}{0.55\hsize}
\begin{center}
\rotatebox{270}{\includegraphics[width=0.75\hsize,clip]
{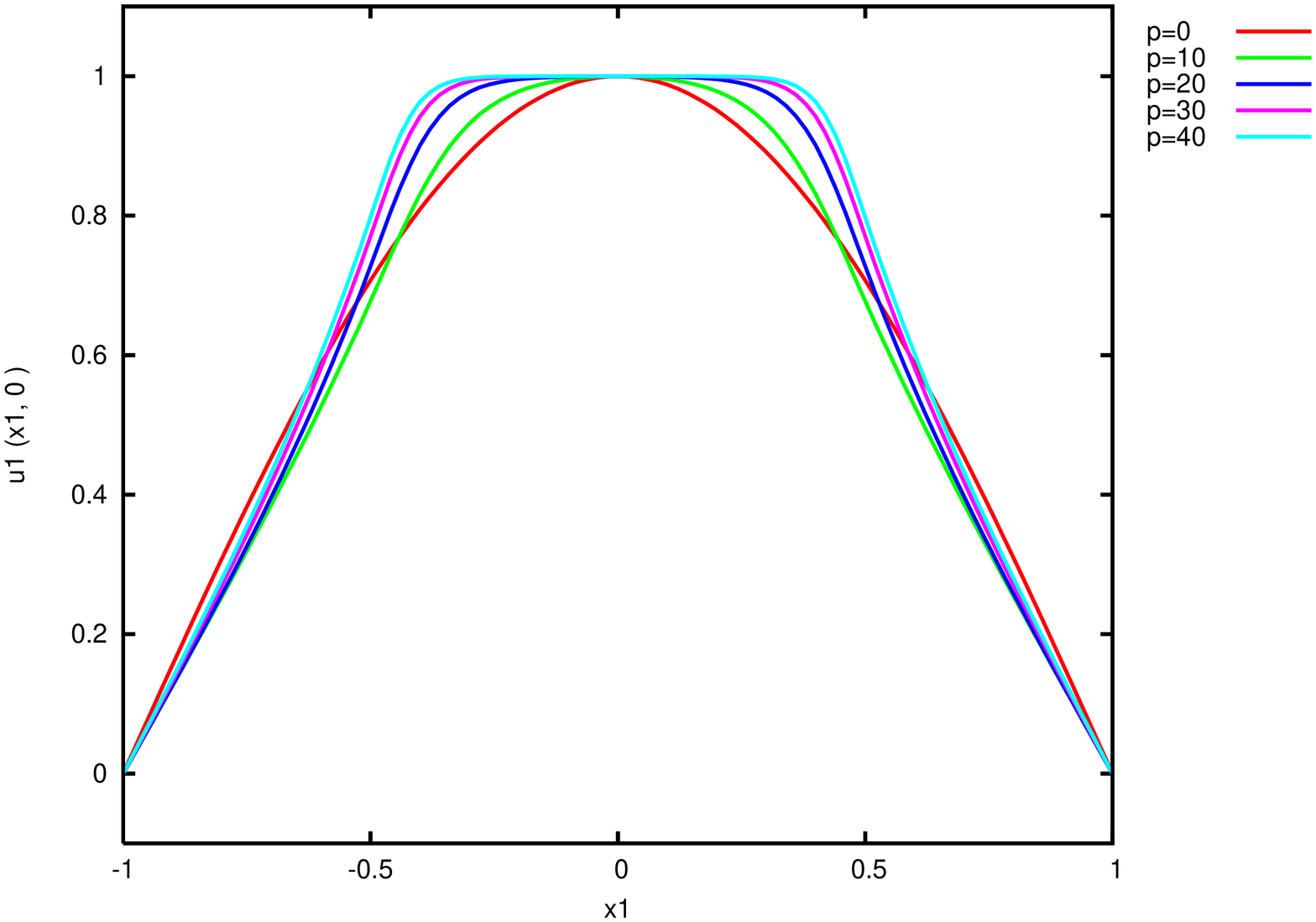}}
\caption{Eigenfunctions of Example~\ref{Sim1} on the line $x_2=0$ for
$p=0,\, 10,\, 20,\, 30,\, 40$.}
\label{Sim1-eigenfunctions}
\end{center}
\end{minipage}
\end{center}
\end{figure}

\begin{Ex}\label{Sim1}{\rm
We set $\ba(\bx):=\balpha (\bx; 1/2)$ and $u_0\equiv 1$.
For $p = 0,\, 10,\, 20,\, 30,\, 40$,
we numerically computed (\ref{parabolic}) 
in the time interval $0\leq t\leq 1.0$
and the results are shown in  
Figures~\ref{Sim1-1d-p00}--\ref{Sim1-eigenfunctions}.
For each $p$, time evolution of $u (\bx,t)$
and a profile of principal eigenfunction $u_1(\bx,p)$
are shown in the figures.

These finite element simulations were computed with time increment
$\tau=1/2000$ on a unstructured triangular mesh of $\Omega$
which consists of 11292 total nodal points (degree of freedom
including boundary points) and 
22182 triangular elements with the average edge length $0.021$.
For the mesh generation, we used FreeFem++ \cite{FFEM}. 

Time evolutions of $u (\bx,t)$ are drawn 
with fixed $x_2=0$,
where the horizontal axis stands for $-1\leq x_1\leq 1$.
In each simulation, for $t\geq 0.2$,
$u (\bx,t)$ decreases monotonically
and keeps almost same profile, i.e.,
$u (\bx,t)$ exhibits the asymptotical form (\ref{asymptotics}).
Comparing these figures, we can see that the decay speed 
of $u (\bx,t)$ becomes extremely slow in cases of $p=20,\, 30,\, 40$.
Especially in case of $p=40$ (Figure~\ref{Sim1-1d-p40}), 
it seems to be almost stationary after $t=0.2$ in the figure.
This is a typical example of 
the exponential decay phenomenon.

The principal eigenfunctions
drawn in the figures are normalized profiles of
$u (\bx,t)$ at $t=1.0$, which are considered as 
profiles of principal eigenfunctions $u_1(\bx,p)$.
The sections of normalized $u_1(\bx,p)$ on the line $x_2=0$ 
are drawn in Figure~\ref{Sim1-eigenfunctions}.
We remark that, in case of $p=0$, 
the principal eigenvalue and eigenfunction
are given by
\[
\lambda_1(0)=\frac{\pi^2}{2},
~~~~~
u_1(\bx,0)=\cos \frac{\pi x_1}{2}\cos \frac{\pi x_2}{2}.
\]
}\end{Ex}

\begin{figure}[htbp]
\begin{center}
\begin{minipage}{0.45\hsize}
\begin{center}
\rotatebox{270}{\includegraphics[width=0.75\hsize,clip]{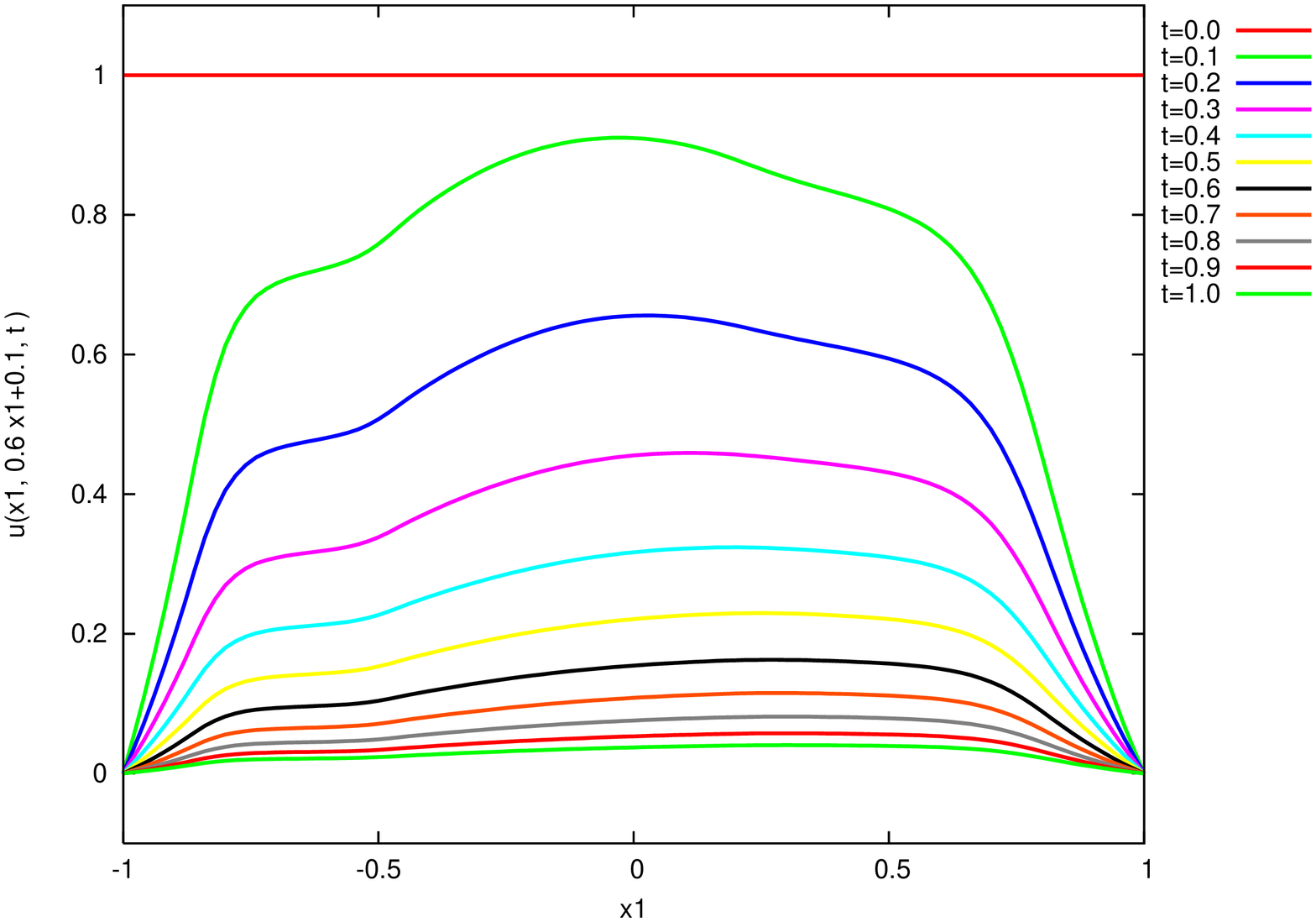}}
\caption{Time evolution of Example~\ref{Sim3} on a line for $p=10$.}
\label{Sim3-1d-p10}
\end{center}
\end{minipage}
\hspace{0.08\hsize}
\begin{minipage}{0.45\hsize}
\begin{center}
\rotatebox{270}{\includegraphics[width=0.75\hsize,clip]{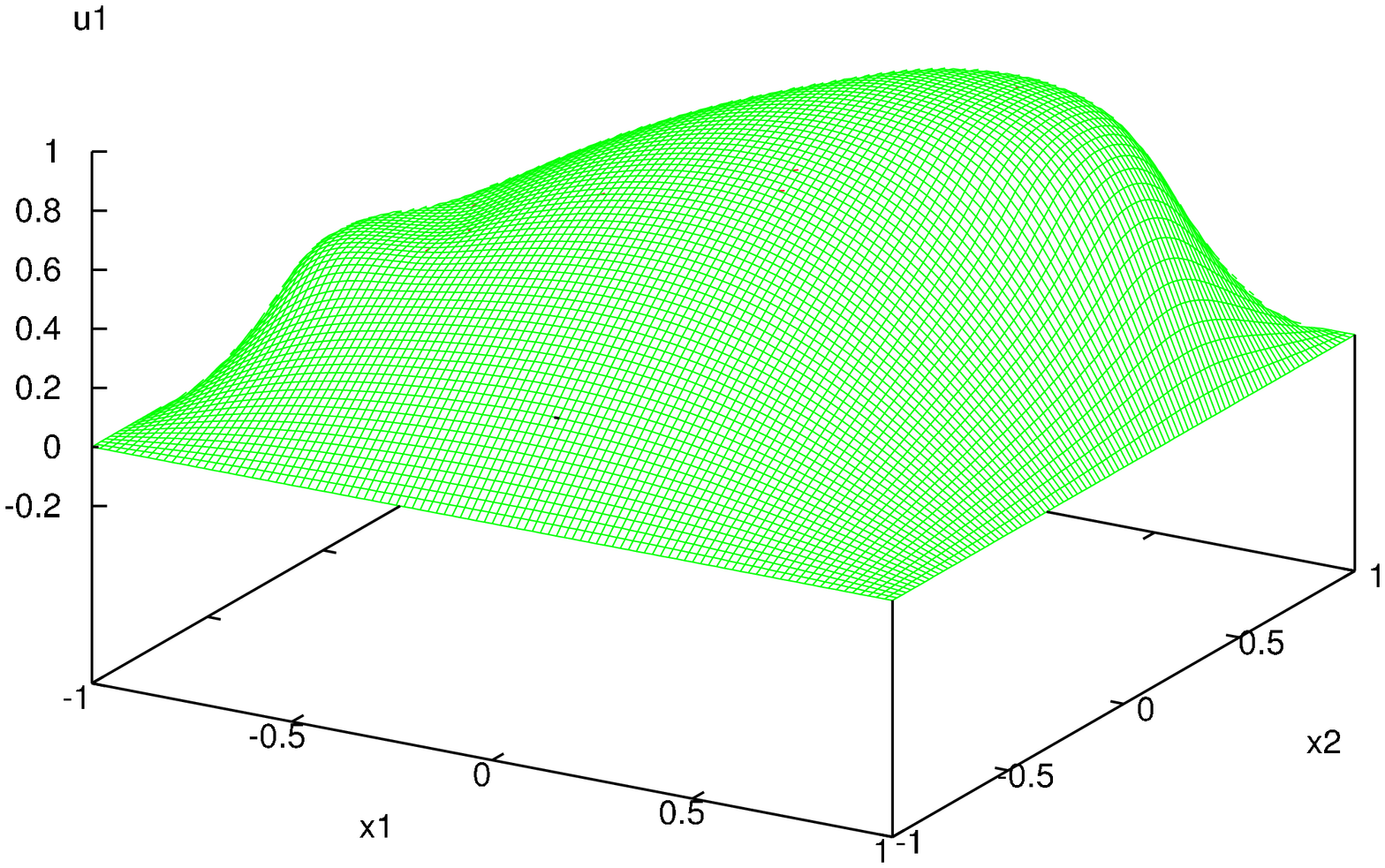}}
\caption{Eigenfunction of Example~\ref{Sim3} for $p=10$.}
\label{Sim3-profile-p10}
\end{center}
\end{minipage}
\end{center}
\begin{center}
\begin{minipage}{0.45\hsize}
\begin{center}
\rotatebox{270}{\includegraphics[width=0.75\hsize,clip]{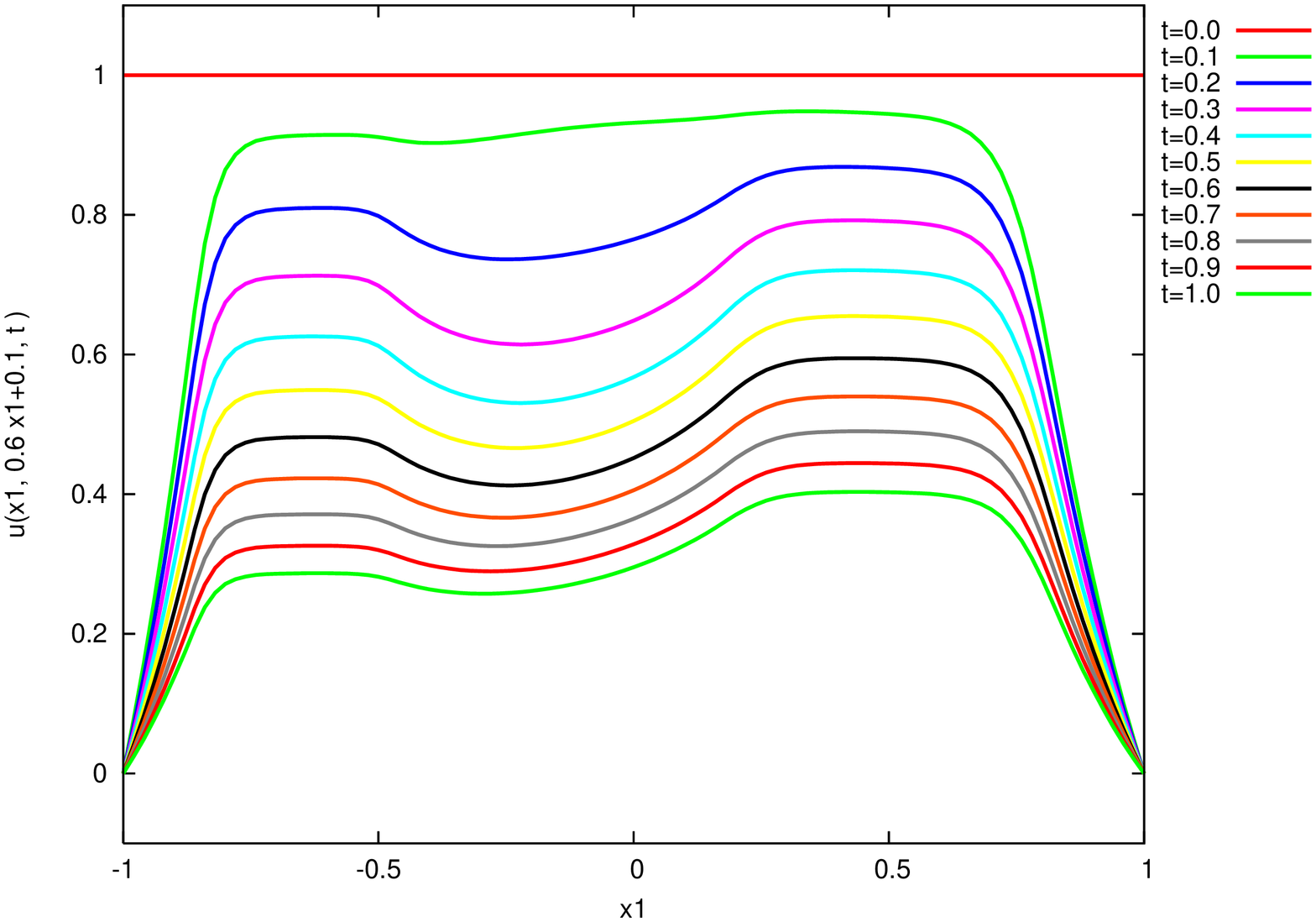}}
\caption{Time evolution of Example~\ref{Sim3} on a line for $p=20$.}
\label{Sim3-1d-p20}
\end{center}
\end{minipage}
\hspace{0.08\hsize}
\begin{minipage}{0.45\hsize}
\begin{center}
\rotatebox{270}{\includegraphics[width=0.75\hsize,clip]{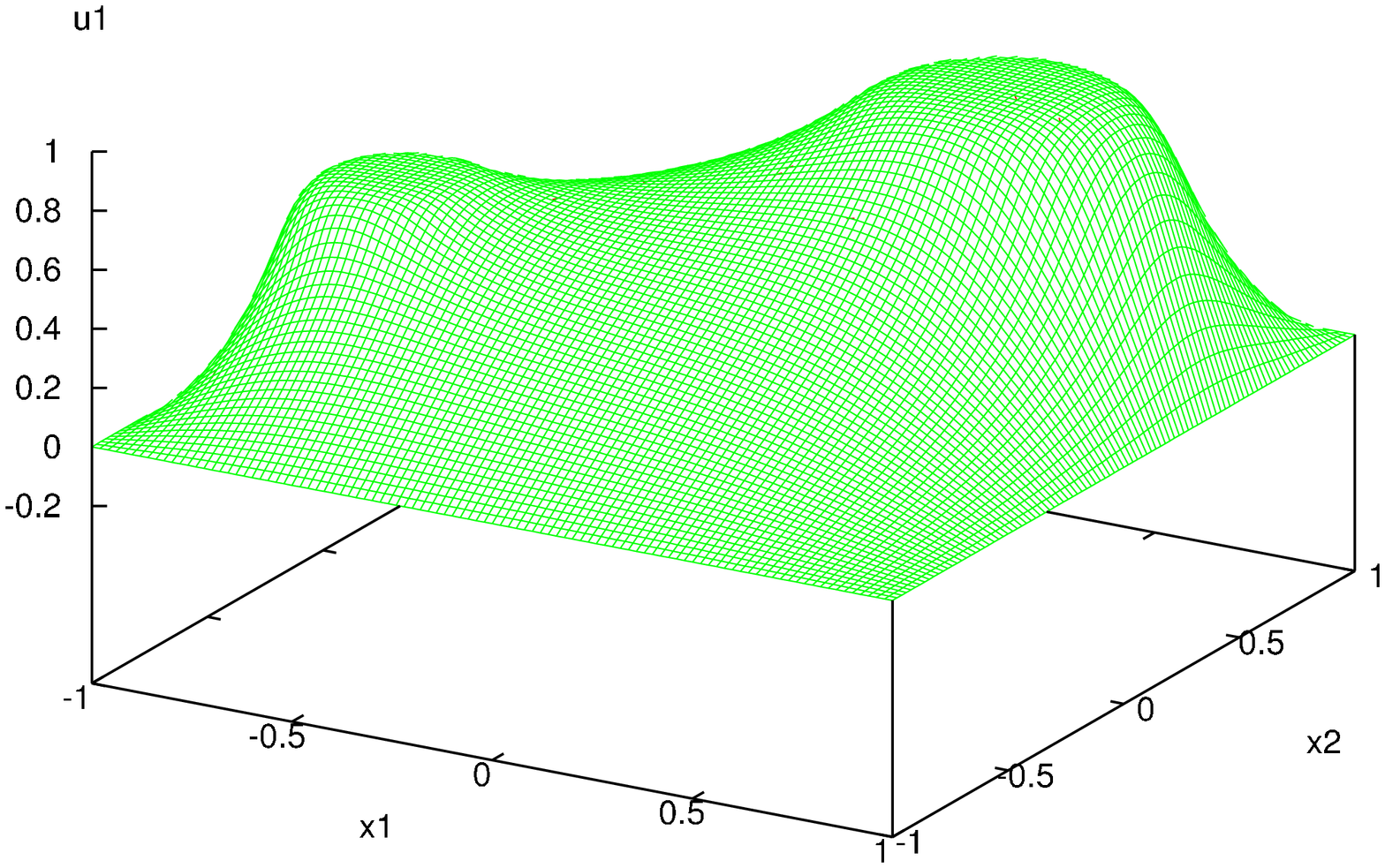}}
\caption{Eigenfunction of Example~\ref{Sim3} for $p=20$.}
\label{Sim3-profile-p20}
\end{center}
\end{minipage}
\end{center}
\begin{center}
\begin{minipage}{0.45\hsize}
\begin{center}
\rotatebox{270}{\includegraphics[width=0.75\hsize,clip]{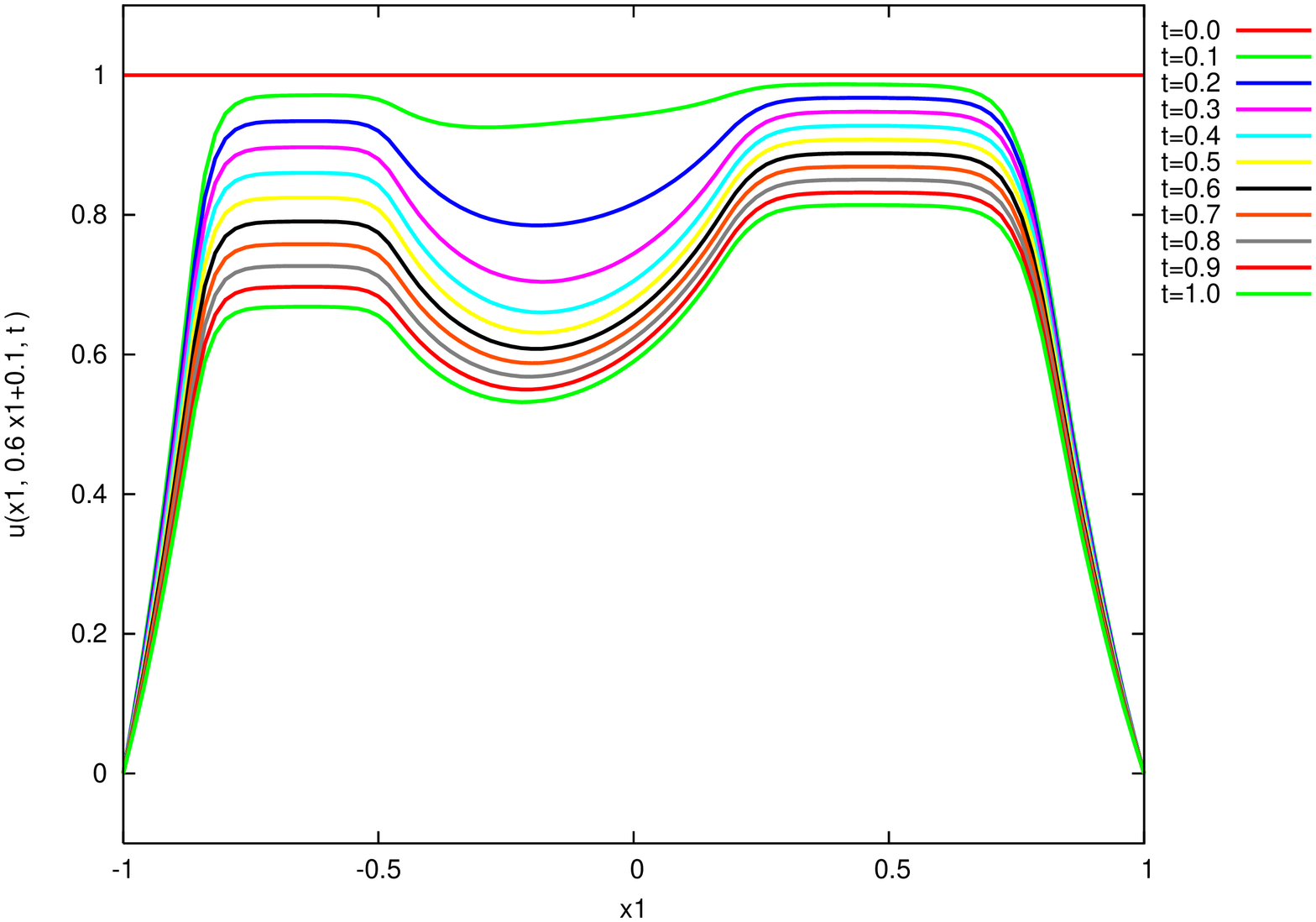}}
\caption{Time evolution of Example~\ref{Sim3} on a line for $p=30$.}
\label{Sim3-1d-p30}
\end{center}
\end{minipage}
\hspace{0.08\hsize}
\begin{minipage}{0.45\hsize}
\begin{center}
\rotatebox{270}{\includegraphics[width=0.75\hsize,clip]{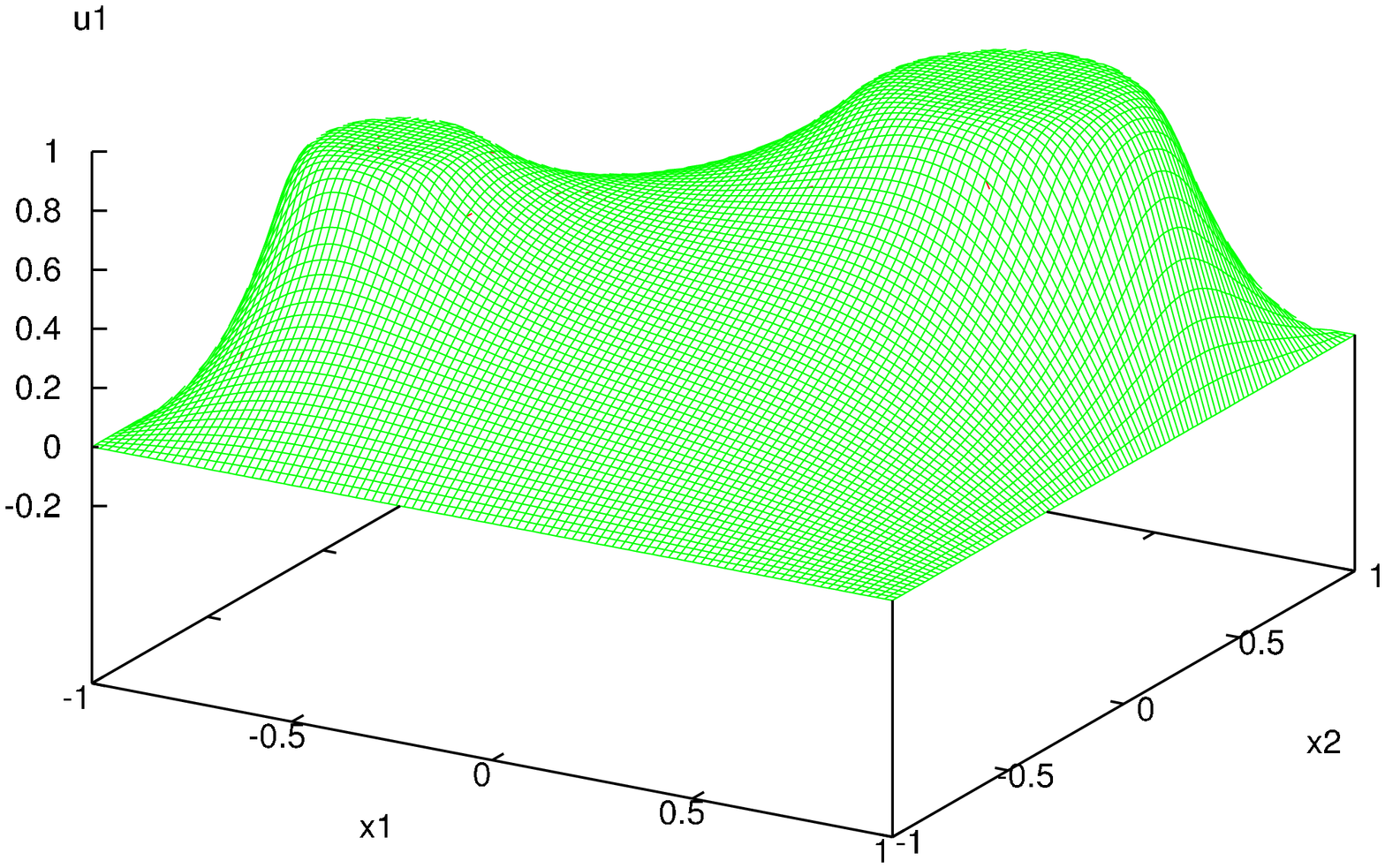}}
\caption{Eigenfunction of Example~\ref{Sim3} for $p=30$.}
\label{Sim3-profile-p30}
\end{center}
\end{minipage}
\end{center}
\end{figure}

\begin{figure}[htbp]
\begin{center}
\begin{minipage}{0.45\hsize}
\begin{center}
\rotatebox{270}{\includegraphics[width=0.75\hsize,clip]{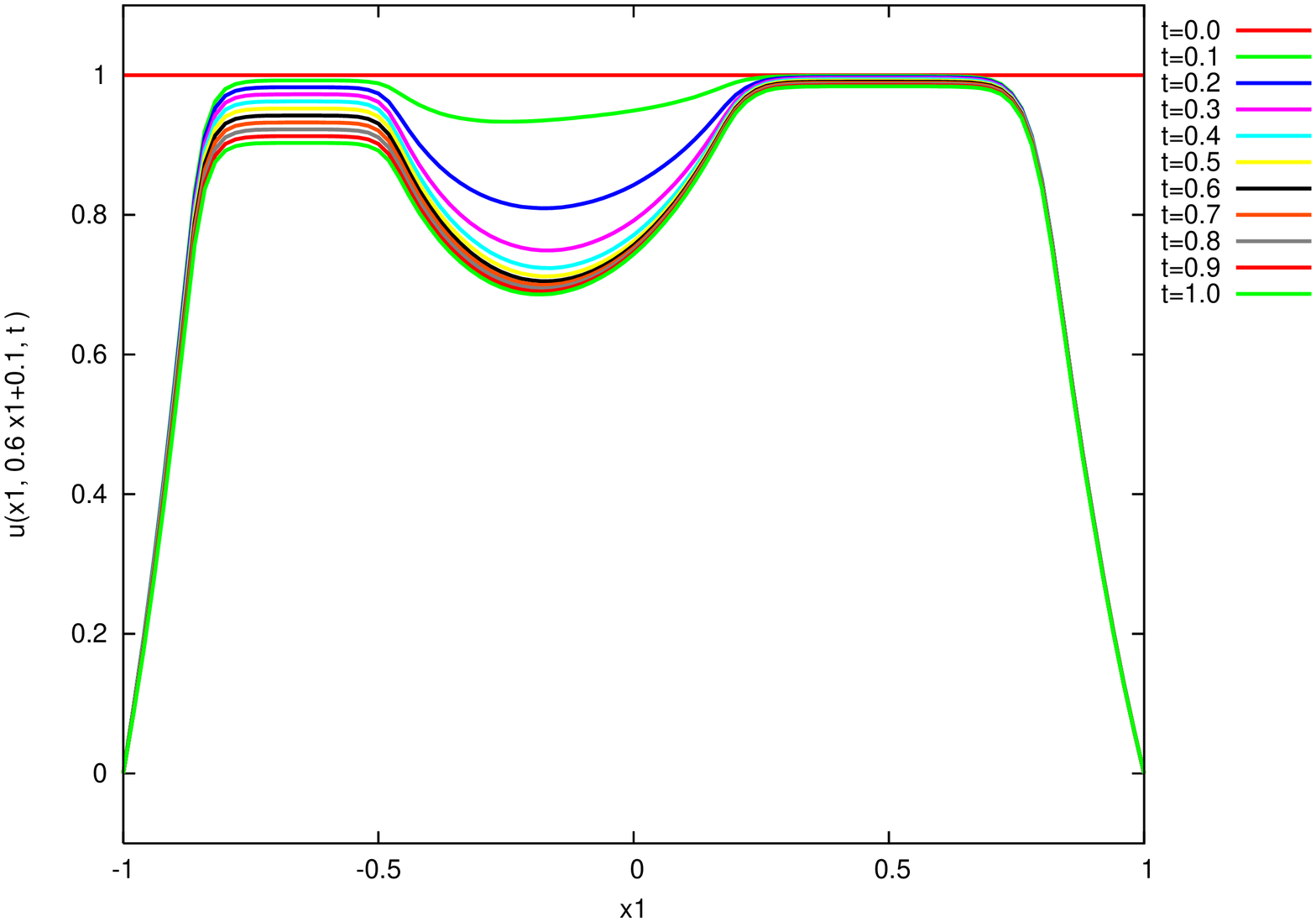}}
\caption{Time evolution of Example~\ref{Sim3} on a line for $p=50$.}
\label{Sim3-1d-p50}
\end{center}
\end{minipage}
\hspace{0.08\hsize}
\begin{minipage}{0.45\hsize}
\begin{center}
\rotatebox{270}{\includegraphics[width=0.75\hsize,clip]{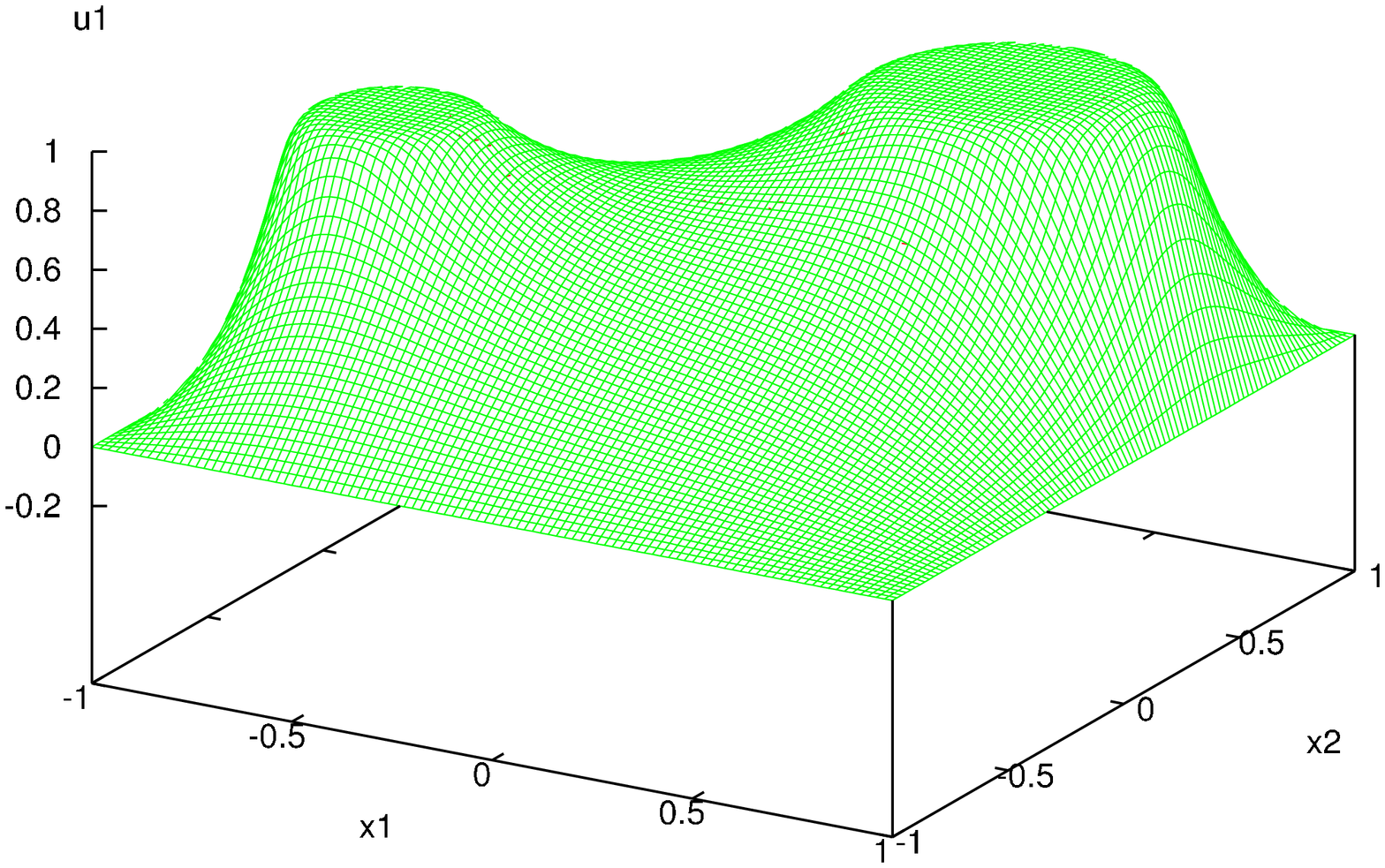}}
\caption{Eigenfunction of Example~\ref{Sim3} for $p=50$.}
\label{Sim3-profile-p50}
\end{center}
\end{minipage}
\end{center}
\begin{center}
\begin{minipage}{0.45\hsize}
\begin{center}
\rotatebox{270}{\includegraphics[width=0.75\hsize,clip]{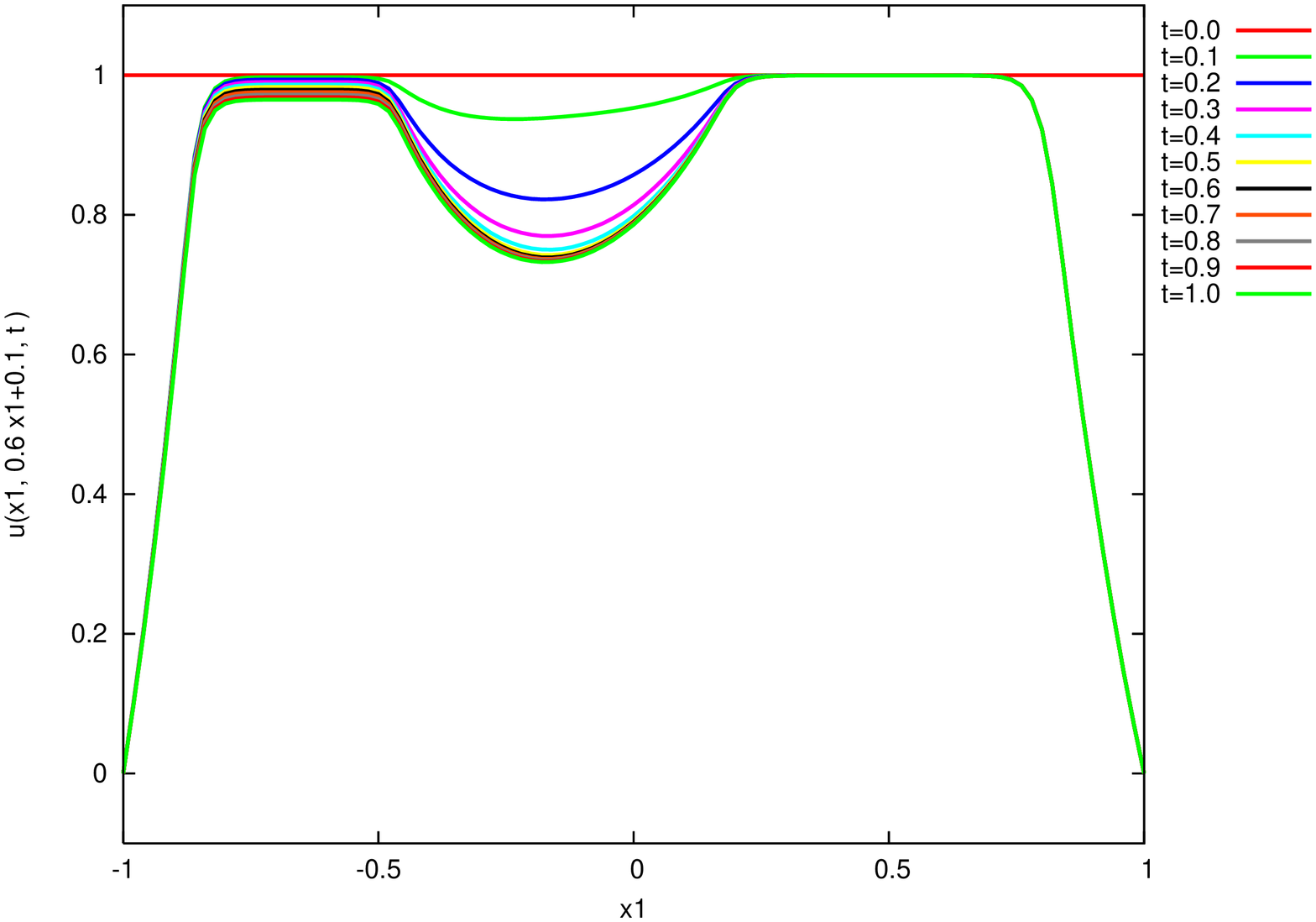}}
\caption{Time evolution of Example~\ref{Sim3} on a line for $p=100$.}
\label{Sim3-1d-p100}
\end{center}
\end{minipage}
\hspace{0.08\hsize}
\begin{minipage}{0.45\hsize}
\begin{center}
\rotatebox{270}{\includegraphics[width=0.75\hsize,clip]{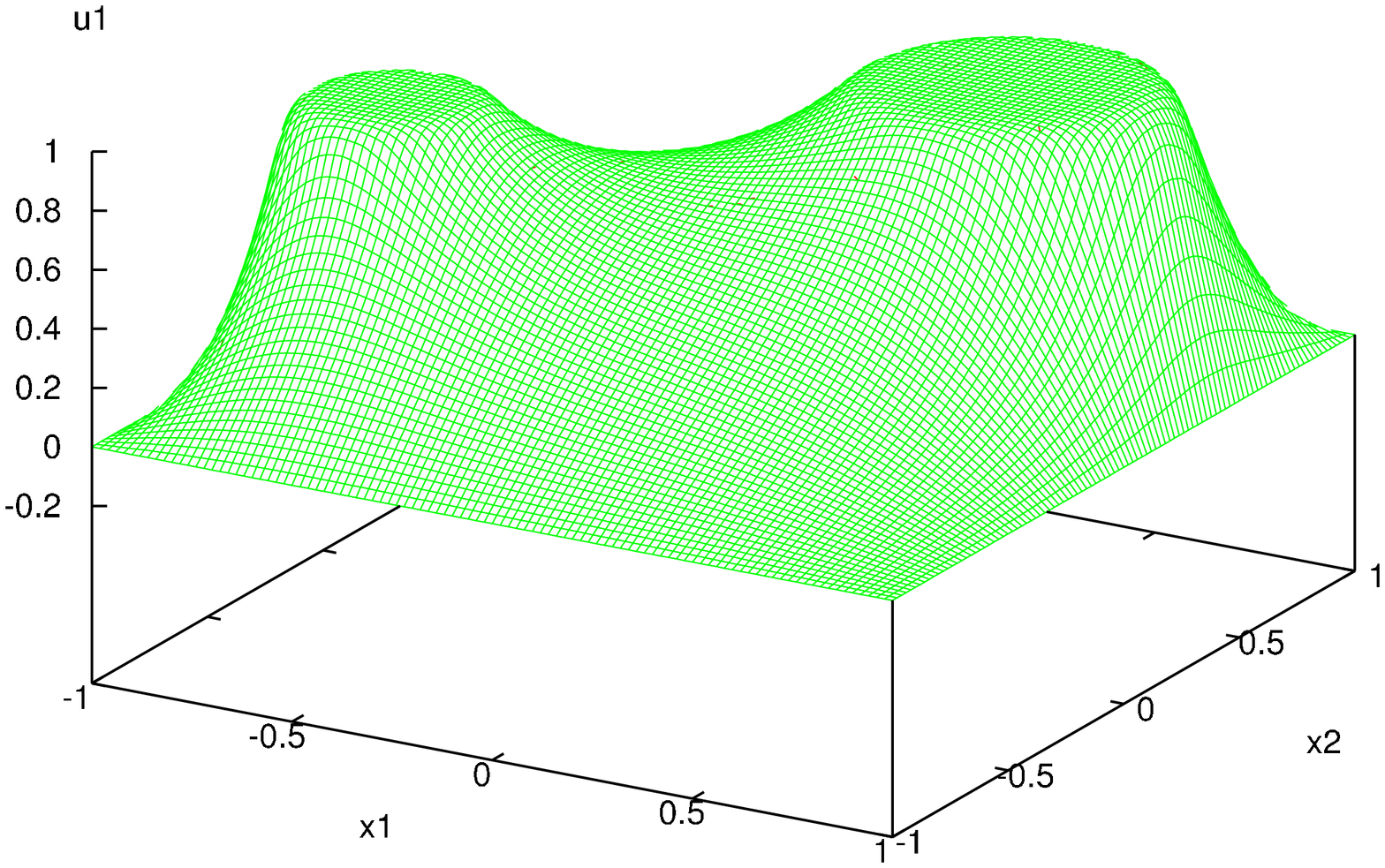}}
\caption{Eigenfunction of Example~\ref{Sim3} for $p=100$.}
\label{Sim3-profile-p100}
\end{center}
\end{minipage}
\end{center}
\begin{center}
\begin{minipage}{0.55\hsize}
\begin{center}
\rotatebox{270}{\includegraphics[width=0.75\hsize,clip]
{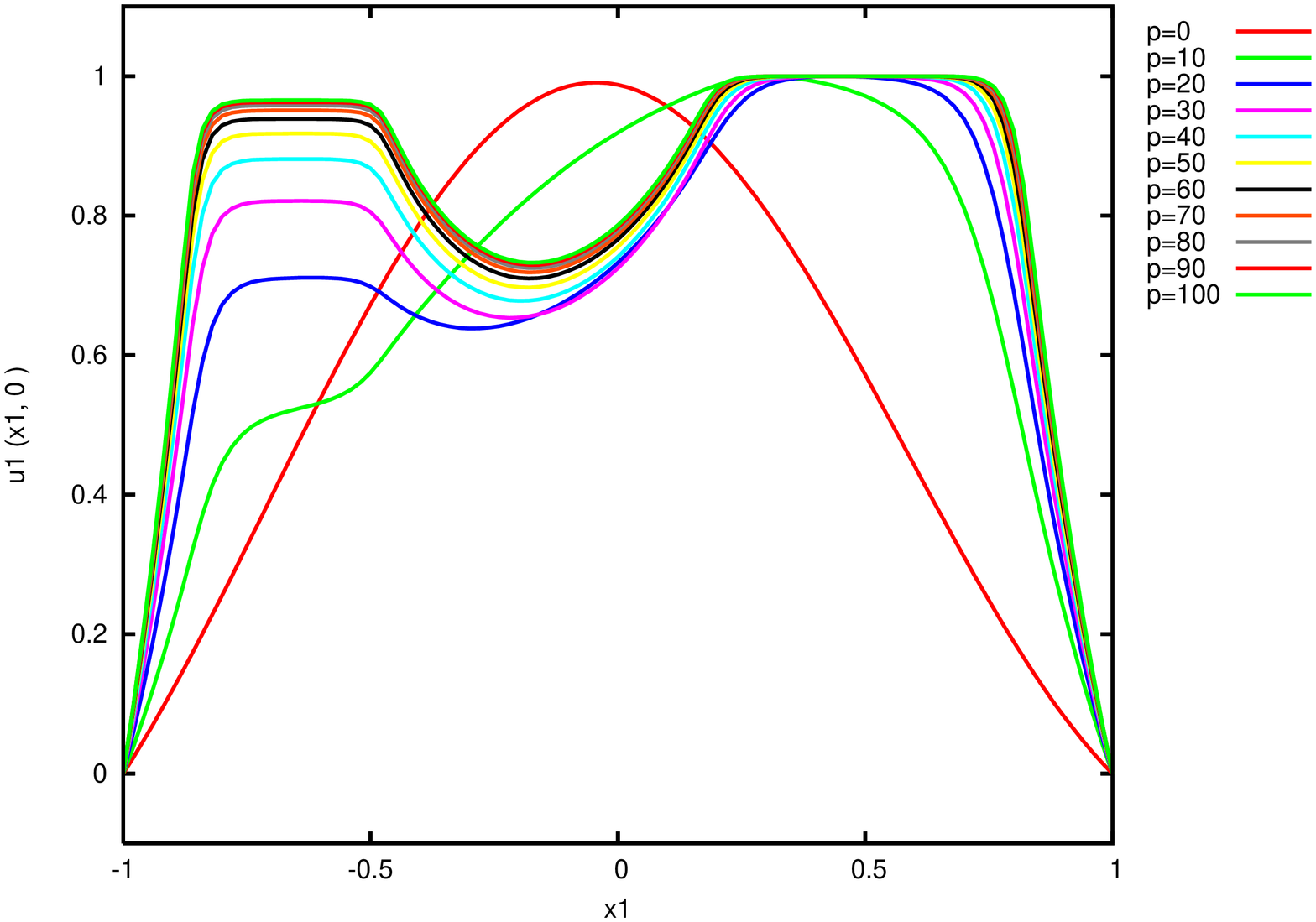}}
\caption{Eigenfunctions of Example~\ref{Sim3} on a line for
$p=0,\, 10,\, 20,\cdots,\, 100$.}
\label{Sim3-eigenfunctions}
\end{center}
\end{minipage}
\end{center}
\end{figure}

\begin{Ex}\label{Sim3}{\rm
We set $\bx_1 := (1/2,\, 2/5)$ and $\bx_2 :=(-2/3,\, -3/10)$
and define
\[
\ba(\bx):=\balpha ( \bx-\bx_1;\,2/5) + 2\balpha ( \bx-\bx_2;\,1/4).
\]
The support of the flow field consists of two disjoint sets
$\{|\bx-\bx_1|\leq 2/5\}$ and  $\{|\bx-\bx_2|\leq 1/4\}$.

For $p = 0,\, 10,\, 20,\cdots,\, 100$,
we numerically computed (\ref{parabolic}) 
in the time interval $0\leq t\leq 1.0$
and the results are shown in  
Figures~\ref{Sim3-1d-p10}--\ref{Sim3-eigenfunctions}.
For each $p=10,\,20,\,30,\,50,\,100$, time evolution of $u (\bx,t)$
and a profile of principal eigenfunction $u_1(\bx,p)$
are shown in the figures.
These finite element simulations were computed with same time increment
and triangulation of $\Omega$ as Example~\ref{Sim1}.

Time evolutions of $u (\bx,t)$ are drawn 
on the line 
$6x_1-10x_2+1=0$ through $\bx_1$ and $\bx_2$,
where the horizontal axis stands for $-1\leq x_1\leq 1$.
Similarly to the previous example, this also exhibits  
the exponential decay phenomenon.

The principal eigenfunctions
drawn in the figures are normalized profiles of
$u (\bx,t)$ at $t=1.0$, which are considered as 
profiles of principal eigenfunctions $u_1(\bx,p)$.
The sections of normalized $u_1(\bx,p)$ on the line $6x_1-10x_2+1=0$
are drawn in Figure~\ref{Sim3-eigenfunctions}.

From this and previous examples, the profiles of the eigenfunction $u_1(\bx,p)$
seems to converge to a limit profile as $p\to\infty$.
In particular, the limit profile seems to be flat on the support of $\ba$
in both cases.
}\end{Ex}
\newpage
\clearpage

\clearpage
\clearpage

\clearpage

\subsection{\large Fundamental tools for eigenvalue problems}
\label{preliminaries} 
\setcounter{equation}{0}

In this section, we collect several fundamental facts and tools
for the eigenvalue problem.
For simple notation, we fix $p=1$ without loss of generality
throughout this section.
\begin{Th}\label{fth}
In a bounded Lipschitz domain $\Omega\subset \R^n$ $(n\geq 1)$,
we assume $\ba\in L^\infty (\Omega,\R^n)$.
Then,
there uniquely exists a positive number $\lambda_1$ and 
$u_1\in H^2_\loc (\Omega)\cap H^1_0(\Omega)$
such that $\lambda=\lambda_1$ and $u=u_1$ satisfy
{\rm (\ref{evp})} with $p=1$ and {\rm (\ref{ucondition})}.
Moreover, $\lambda_1$ is given by
the min-max formula:
\[
\lambda_1=\max_{\varphi>0}\,
\inf_{\sbx\in\Omega}
\left(
\frac{-\Delta\varphi(\bx)+\ba(\bx)\cdot\nabla \varphi(\bx)}
{\varphi(\bx)}
\right),
\]
where $\max_{\varphi>0}$ is taken over all
positive $\varphi\in W^{2,n}_\loc (\Omega)$.
\end{Th}
We remark that $\inf$ and $\sup$ always denote essential infimum
and supremum in this paper.
For the proof of the theorem, see 
Section~2.8 of \cite{P-W84} and \cite{B-N-V94}.
Since
$\lambda_1=\min_{\lambda} \mbox{Re}\lambda$ holds
for any complex eigenvalue $\lambda$ of $-\Delta +\ba\cdot\nabla$,
$\lambda_1$ and $u_1$ are called the principal 
eigenvalue and eigenfunction.
The next corollary immediately follows from Theorem~\ref{fth}.
\begin{Cor}\label{fcor}
Under the condition of Theorem~\ref{fth}, we consider
a Lipschitz subdomain $\Omega'\subset \Omega$ and 
define $\lambda_1'>0$ as the principal eigenvalue for
$-\Delta u+\ba\cdot \nabla u=\lambda_1'u$ in $\Omega'$
with the zero Dirichlet boundary condition on $\partial\Omega'$.
Then $\lambda_1\leq \lambda_1'$ holds.
\end{Cor}

Throughout the following sections, we assume the existence of a
velocity potential of $\ba$:
\begin{equation}\label{assumeb}
~^\exists b\in W^{2,\infty}(\Omega)
~\mbox{s.t.}~
\ba(\bx)=\nabla b(\bx)
~~~~~
(\bx\in\Omega).
\end{equation}
Since the existence of a
velocity potential is not essential for the exponential decay phenomena
which is our aim in this paper (see \cite{Fri73}, for example),
the assumption (\ref{assumeb}) is rather technical
but allows us energy-based arguments with the help of
the following Liouville transform. 

Under the condition (\ref{assumeb}) and for $q\in L^\infty(\Omega)$,
we consider the following three elliptic eigenvalue
problems with the Dirichlet boundary condition.
\begin{eqnarray}
&&\left\{
\begin{array}{ll}
-\Delta u (\bx)+\ba(\bx)\cdot\nabla u (\bx)=\lambda u (\bx)
&(\bx\in\Omega)\\
~&\\
u (\bx)=0
&(\bx\in\partial\Omega)\\
\end{array}
\right. \label{eigen-v}\\~\nonumber \\
&&\left\{
\begin{array}{ll}
-\Delta v (\bx)-\div (v (\bx)\ba(\bx))=\lambda v (\bx)
&(\bx\in\Omega)\\
~&\\
v (\bx)=0
&(\bx\in\partial\Omega)\\
\end{array}
\right. \label{eigen-u}\\~\nonumber \\
&&\left\{
\begin{array}{ll}
-\Delta w(\bx)+q(\bx) w(\bx)=\lambda w(\bx)
&(\bx\in\Omega)\\
~&\\
w(\bx)=0
&(\bx\in\partial\Omega)\\
\end{array}
\right. \label{eigen-w}
\end{eqnarray}

\begin{Prop}[Liouville transform]\label{equivalence}
Suppose the condition {\rm (\ref{assumeb})}.
Then
the above three eigenvalue
problems are equivalent to each other under the
following relations:
\[
e^{\frac{b}{2}}v =e^{-\frac{b}{2}}u=w,
~~~~~
q= -\frac{1}{2}\Delta b+\frac{1}{4}|\nabla b|^2.
\]
\end{Prop}
\noindent
We omit the proof, since it is shown by straightforward
substitution.
We remark that the same Liouville transform is valid
even for the time dependent problem (\ref{parabolic}).

For the self-adjoint eigenvalue problem (\ref{eigen-w})
with $q\in L^\infty(\Omega)$,
there is a well-known min-max principle for the Rayleigh quotient.
We define
\begin{equation}\label{Lq}
L_q:=- \Delta +q,~~~~~
\mbox{\rm Dom}(L_q)=H^1_0(\Omega),
\end{equation}
\[
J_q(w):=
\frac{\mbox{}_{H^{-1}(\Omega)}\langle
L_qw,w\rangle_{H^1_0(\Omega)}}
{(w,w)_{L^2(\Omega)}}
=
\frac{\int_\Omega ( |\nabla w(\bx)|^2+q(\bx)w(\bx)^2)d\bx}
{\int_\Omega w(\bx)^2d\bx}
~~~~~
(w\in H^1_0(\Omega)).
\]

\begin{Th}\label{ONS}
For $q\in L^\infty(\Omega)$, there exists
a complete orthonormal system of $L^2(\Omega)$,
$\{w_k\}_{k\in\N}\subset H^1_0(\Omega)\cap H^2_\loc (\Omega)$
and $\{\lambda_k\}_{k\in\N}\subset \R$
$(\lambda_1<\lambda_2\leq\cdots)$ such that
\[
L_qw_k=\lambda_kw_k
~~~~~
(k\in\N).
\]
The $k$th eigenvalue $\lambda_k$ is characterized by the
min-max formula:
\[
\lambda_k=\min_{\dim X=k}
\max_{w\in X\setminus \{0\}}
J_q(w),
\]
where the minimum is taken over whole 
$k$-dimensional subspace $X$ in $H^1_0(\Omega)$.
In particular, the principal
eigenvalue $\lambda_1$ is simple and given by
\[
\lambda_1=\min_{w\in H^1_0(\Omega)\setminus \{0\}}
J_q(w).
\]
Moreover, an eigenfunction $w\in H^1_0(\Omega)\cap H^2_\loc (\Omega)$
corresponds to the principal eigenvalue if and only if
$w$ has no sign change in $\Omega$.
\end{Th}
\noindent
For the proof of this theorem, 
see  Section~8.12 of \cite{G-T83}
and \cite{B-B-B-S06} etc.

From the last assertion of this theorem,
without loss of generality, 
we can assume that
\[
w_1(\bx)>0~~~(\bx\in\Omega).
\]

\begin{Th}[comparison theorem for principal eigenvalues]\label{cth}
For $q,~\tilde{q}\in L^\infty(\Omega)$,
the principal eigenvalues of $L_q$ and $L_{\tilde{q}}$
are denoted by $\lambda_1$ and $\tilde{\lambda}_1$,
respectively.
Then
\[
\einf_{\Omega}(q-\tilde{q})
\leq
\lambda_1-\tilde{\lambda}_1
\leq
\esup_{\Omega}(q-\tilde{q}).
\]
In particular,
if $q(\bx)\geq\tilde{q}(\bx)$ $(\bx\in\Omega)$
and $\| q-\tilde{q}\|_{L^\infty(\Omega)}\neq 0$, then
$\lambda_1>\tilde{\lambda}_1$.
\end{Th}
\noindent
{\it Proof.}
Let $w_1$ and $\tilde{w}_1$ denote
the eigenfunctions corresponding to $\lambda_1$
and $\tilde{\lambda}_1$, respectively.
We have
\begin{eqnarray*}
\lambda_1
&=&
\min_{w\in H^1_0(\Omega),w\not\equiv 0}
J_q(w)
\leq
J_q(\tilde{w}_1)
=
J_{\tilde{q}}(\tilde{w}_1)
+\left(
J_q(\tilde{w}_1)-J_{\tilde{q}}(\tilde{w}_1)
\right)\\
&=&
\tilde{\lambda}_1
+\frac{\int_\Omega (q(\bx)-\tilde{q}(\bx))\tilde{w}_1(\bx)^2d\bx}
{\int_\Omega \tilde{w}_1(\bx)^2d\bx}
\leq
\tilde{\lambda}_1
+\esup_{\Omega}(q-\tilde{q}).
\end{eqnarray*}
In the same way, we also have
\begin{equation}\label{lam*}
\tilde{\lambda}_1
\leq
\lambda_1
+\esup_{\Omega}(\tilde{q}-q)
=
\lambda_1
-\einf_{\Omega}(q-\tilde{q}).
\end{equation}
These inequalities imply the first assertion. Furthermore, if 
$q\geq\tilde{q}$ and $q\not\equiv \tilde{q}$ then,
from the positivity of $w_1$ and $\tilde{w}_1$, we have
\[
\lambda_1(w_1,\tilde{w}_1)_{L^2(\Omega)}
=(L_qw_1,\tilde{w}_1)_{L^2(\Omega)}
=(w_1,L_q\tilde{w}_1)_{L^2(\Omega)}
>(w_1,L_{\tilde{q}}\tilde{w}_1)_{L^2(\Omega)}
=\tilde{\lambda}_1(w_1,\tilde{w}_1)_{L^2(\Omega)},
\]
and $\lambda_1>\tilde{\lambda}_1$ follows.
\qed \\

We define the elliptic operator $K_b$ by
\[
K_b:=- \Delta +\nabla b\cdot\nabla,~~~~~
\mbox{\rm Dom}(K_b)=H^1_0(\Omega),
\]
and denote by $K_b^*$ its adjoint operator
with respect to $L^2(\Omega)$:
\[
K_b^*:=- \Delta -\div (\cdot \nabla b),~~~~~
\mbox{\rm Dom}(K_b^*)=H^1_0(\Omega).
\]
Although the following results concerning
the operator $K_b^*$ will not be used in our
analysis,
we include them below
for our systematic description.
Another reason to include them 
is that the eigenvalue problem
(\ref{eigen-u}) corresponding to $K_b^*$ 
is related to various important applications
as well as (\ref{eigen-v}) for $K_b$.
An example of such applications
will be shown in Section~\ref{chemotaxis}. 

We introduce the following weighted inner product of
$L^2(\Omega)$;
\[
(v ,u)_{L^2_b(\Omega)}:=\int_\Omega
e^{b(\sbx)}v (\bx)u (\bx)d\bx,
~~~~~
(v ,~u\in L^2(\Omega)),
\]
and we define $L^2_b(\Omega):=
(L^2(\Omega),~(\cdot,\cdot)_{L^2_b(\Omega)})$
which denotes a Hilbert space
$L^2(\Omega)$ with this inner product.
We also define
\[
\begin{array}{ll}
\DS{I_b(u):=
\frac{\int_\Omega e^{-b(\sbx)}
|\nabla u (\bx)|^2d\bx}
{\int_\Omega e^{-b(\sbx)}u (\bx)^2d\bx}}
&
(u\in H^1_0(\Omega),~u\not\equiv 0),\\~\\
\DS{I_b^*(v ):=
\frac{\int_\Omega e^{b(\sbx)}
\{ |\nabla v (\bx)|^2
-\Delta b(\bx)v (\bx)^2\}d\bx}
{\int_\Omega e^{b(\sbx)}v (\bx)^2d\bx}}~~~
&
(v \in H^1_0(\Omega),~v \not\equiv 0).
\end{array}
\]
Then we have the following proposition.
\begin{Prop}
We suppose that $u ,\,v,\,\varphi \in H^1_0(\Omega)\cap H^2(\Omega)$
$(u,\,v\not\equiv 0)$.
Then we have
\[
(K_b u,\varphi)_{L^2_{-b}(\Omega)}
=(u,K_b\varphi)_{L^2_{-b}(\Omega)},
~~~~~
I_b(u)=
\frac{(K_bu,u)_{L^2_{-b}(\Omega)}}
{(u,u)_{L^2_{-b}(\Omega)}},
\]
\[
(K_b^*v ,\varphi)_{L^2_b(\Omega)}
=(v ,K_b^*\varphi)_{L^2_b(\Omega)},
~~~~~~
I_b^*(v )=
\frac{(K_b^*v ,v )_{L^2_b(\Omega)}}
{(v ,v )_{L^2_b(\Omega)}}.
\]
Namely, $K_b$ and $K_b^*$ are selfadjoint in
$L^2_{-b}(\Omega)$ and $L^2_b(\Omega)$,
respectively, and their Rayleigh quotients
are given by $I_b(u)$ and $I_b^*(v )$.
\end{Prop}
\noindent
{\it Proof.}
For $\varphi=u$ or $v$, 
by the integration by parts,
we obtain
\[
\int_\Omega e^{\mp b}
|\nabla \varphi|^2d\bx
=
-\int_\Omega
\div \left( e^{\mp b}\nabla \varphi\right)
\varphi \,d\bx
=
\int_\Omega
e^{\mp b}
\left(-\Delta \varphi \pm \nabla b\cdot \nabla \varphi\right)
\varphi \,d\bx .
\]
Hence, we have
\[
I_b(u)=
\frac{
\int_\Omega e^{-b}
(-\Delta u + \nabla b\cdot \nabla u)\,u
d\bx
}
{
(u,u)_{L^2_{-b}(\Omega)}
}
=
\frac{
\int_\Omega e^{-b}(K_bu)\,u d\bx
}
{
(u,u)_{L^2_{-b}(\Omega)}
}
=
\frac{(K_bu,u)_{L^2_{-b}(\Omega)}}
{(u,u)_{L^2_{-b}(\Omega)}},
\]
\[
I_b^*(v )=
\frac{
\int_\Omega e^{b}
\{(-\Delta v  - \nabla b\cdot \nabla v )\,v 
-(\Delta b)\,v^2\}d\bx
}
{
(v ,v )_{L^2_b(\Omega)}
}
=
\frac{
\int_\Omega e^{b}(K_b^*v )\,v  d\bx
}
{
(v ,v )_{L^2_b(\Omega)}
}
=
\frac{(K_b^*v ,v )_{L^2_b(\Omega)}}
{(v ,v )_{L^2_b(\Omega)}}.
\]

The self-adjointness of $K_b$ in $L^2_{-b}(\Omega)$
and the one of $K_b^*$ in $L^2_b(\Omega)$
are also checked by direct calculations.
\qed \\

We have the following characterizations of
the eigenvalues
of (\ref{eigen-v}) and (\ref{eigen-u}).
\begin{Th}\label{ONS2}
For $b\in W^{2,\infty}(\Omega)$, there exist
$\{\lambda_k\}_{k\in\N}\subset \R$
$(\lambda_1<\lambda_2\leq\cdots)$,
and $\{u_k\}_{k\in\N}$, $\{v_k\}_{k\in\N}
\subset H^1_0(\Omega)\cap H^2_\loc (\Omega)$
such that
\[
K_bu_k=\lambda_k u_k,
~~~~~
K_b^*v_k=\lambda_k v_k
~~~~~
(k\in\N),
\]
and that $\{u_k\}_{k\in\N}$ and
$\{v_k\}_{k\in\N}$ are complete orthonormal systems
of $L^2_{-b}(\Omega)$ and $L^2_b(\Omega)$,
respectively.
The $k$th common eigenvalue $\lambda_k$ is characterized by the
min-max formula:
\[
\lambda_k=\min_{\dim X=k}
\max_{u\in X\setminus \{0\}}
I_b(u)=
\min_{\dim X=k}
\max_{v \in X\setminus \{0\}}
I_b^*(v ),
\]
where the minimum is taken over whole 
$k$-dimensional subspace $X$ in $H^1_0(\Omega)$.
In particular, the common principal
eigenvalue $\lambda_1$ is simple 
and given by
\[
\lambda_1
=\min_{u\in H^1_0(\Omega)\setminus \{0\}}I_b(u)
=\min_{v \in H^1_0(\Omega)\setminus \{0\}}I_b^*(v ).
\]
Moreover, an eigenfunction $u$ for $K_b$ $(v $ for $K_b^*)$
corresponds to the principal eigenvalue if and only if
$u$ $(v )$ has no sign change in $\Omega$.
\end{Th}
\noindent
{\it Proof.}
The assertions follow from
Proposition~\ref{equivalence}
and Theorem~\ref{ONS}.
\qed \\

From the last assertion of this theorem,
without loss of generality, 
we can assume that
\[
u_1(\bx)>0,~~~~~v_1(\bx)>0~~~(\bx\in\Omega).
\]
As a last remark in this section,
we introduce the following proposition.
We also omit the proof since it is shown by a direct calculation.
\begin{Prop}\label{rect}
Let $\Omega:=\Pi_{i=1}^n \cI_i\subset \R^n$,
where $\cI_i$ is a bounded open interval for $i=1,\cdots,n$,
and let $\ba (\bx)=~(a_1(x_1),a_2(x_2),\cdots,a_n(x_n))^\rT$
for $\bx=(x_1,\cdots,x_n)^\rT \in \Omega$
with $a_i\in L^\infty(\cI_i)$.
Then the principal eigenvalue $\lambda_1$ and the eigenfunction
$u_1(\bx)$ of {\rm (\ref{eigen-v})} are given by
the following formula:
\[
\lambda_1=\sum_{i=1}^n \lambda_1^{(i)},
~~~~~
u_1(\bx):=\prod _{i=1}^n \varphi_i(x_i)
~~~
(\bx\in \Omega),
\]
where
$\lambda_1^{(i)}$ and $\varphi_i$ are defined by
\[
\left\{
\begin{array}{ll}
-\varphi_i''(x)+a_i(x)\varphi_i'(x)=\lambda_1^{(i)}\varphi_i(x)
&(x\in I_i)\\~\\
\varphi_i(x)=0
&(x\in\partial I_i)\\~\\
\varphi_i(x)>0&(x\in I_i).
\end{array}
\right.
\]
\end{Prop}

\subsection{\large Asymptotic behaviour of principal eigenvalues}
\label{asymptotic}
\setcounter{equation}{0}

In this section, we consider the singular perturbation 
problem of the principal
eigenvalues (\ref{evp})
under the velocity potential condition 
(\ref{assumeb}).
Henceforth we assume that $p$ is a positive parameter.

Due to Proposition~\ref{equivalence},
(\ref{evp}) is equivalent to
each of the following three 
eigenvalue problems under the zero Dirichlet boundary condition:
\[
\begin{array}{lll}
K_{pb}u=-\Delta u+p\nabla b\cdot\nabla u=\lambda u
&\mbox{in}&\Omega,\\
~\\
K_{pb}^*v =-\Delta v -p\,\div (v \nabla b)=\lambda v 
&\mbox{in}&\Omega,\\
~\\
L_{q(p)}w=-\Delta w+q(p) w=\lambda w
&\mbox{in}&\Omega,
\end{array}
\]
where we define
\begin{equation}\label{qxp}
q(\bx,p):=-\frac{p}{2} \div \ba(\bx)+\frac{p^2}{4}|\ba(\bx)|^2,
\end{equation}
and we abbreviate $q(\cdot,p)$($\in L^\infty(\Omega)$) as 
$q(p)$.

Under the condition (\ref{assumeb}),
the $k$th eigenvalue (their multiplicities are counted)
of these equivalent eigenvalue problems 
is denoted by $\lambda_k(p)$
for a parameter $p>0$.

We start from the following simple consequence of
the comparison theorem.

\begin{Th}
Let $\lambda_\Omega>0$ be the principal eigenvalue 
of $-\Delta_D$, which is Laplacian 
with the zero Dirichlet boundary 
condition.
Under the condition {\rm (\ref{assumeb})}, following two estimates hold:
\[
\lambda_1(p)\geq \lambda_\Omega
-\frac{p}{2}
\sup_{\sbx\in\Omega} (\div \ba(\bx)) 
~~~~~(p>0).
\]
\[
\frac{1}{4} \inf_{\sbx\in\Omega}|\ba(\bx)|^2
\leq \,\liminf_{p\to\infty} \,\frac{\lambda_1(p)}{p^2}\,
\leq \,\limsup_{p\to\infty} \,\frac{\lambda_1(p)}{p^2}\,
\leq
\frac{1}{4} \sup_{\sbx\in\Omega}|\ba(\bx)|^2.
\]
\end{Th}
\noindent
{\it Proof.}
Applying Theorem~\ref{cth} to $L_{q(p)}$
and $L_0=-\Delta_D$, we have
\begin{equation}\label{compare1}
\inf_{\sbx\in\Omega}q(\bx,p)
\leq
\lambda_1(p)- \lambda_\Omega
\leq
\sup_{\sbx\in\Omega}q(\bx,p),
\end{equation}
and
\begin{equation}\label{infest}
\lambda_1(p)
\geq
\lambda_\Omega
+\inf
\left(
-\frac{p}{2} \div \ba+\frac{p^2}{4}|\ba|^2
\right)\geq 
\lambda_\Omega
-\frac{p}{2} \sup (\div \ba)+\frac{p^2}{4}\inf |\ba|^2,
\end{equation}
\begin{equation}\label{supest}
\lambda_1(p)
\leq
\lambda_\Omega
+\sup
\left(
-\frac{p}{2} \div \ba+\frac{p^2}{4}|\ba|^2
\right)
\leq
\lambda_\Omega
-\frac{p}{2} \inf (\div \ba)+\frac{p^2}{4}\sup |\ba|^2.
\end{equation}
The first assertion of the theorem follows from 
(\ref{infest}). Dividing (\ref{infest}) and (\ref{supest}) by $p^2$
and taking the limit $p\to\infty$, 
we can also derive the second assertion.
\qed \\

\begin{Cor}
There is no exponential decay phenomenon if 
$\inf|\ba(\bx)|>0$
or 
$\sup \div \ba(\bx)\leq 0$.
\end{Cor}

The estimate (\ref{compare1}) can be improved as follows.
\begin{Lem}\label{lemp12}
For $p_1>p_2\geq 0$, we have
\begin{equation}\label{compare2}
\einf_{\bx\in\Omega}
q(\bx,p_1+p_2)
\leq 
\frac{p_1+p_2}{p_1-p_2}(\lambda_1(p_1)-\lambda_1(p_2))
\leq
\esup_{\bx\in\Omega}
q(\bx,p_1+p_2).
\end{equation}
\end{Lem}
\noindent
{\it Proof.}
From the definition of $q(\bx,p)$ (\ref{qxp}), 
the equality
\[
q(\bx,p_1)-q(\bx,p_2)
=\frac{p_1-p_2}{p_1+p_2}q(\bx,p_1+p_2),
\]
holds. Hence, the assertion follows from
Theorem~\ref{cth}.
\qed \\

From (\ref{compare1}) and (\ref{compare2}),
it is natural to consider the following
condition for $\ba$:
\begin{equation}\label{condition1}
^\exists p_0>0
~\mbox{s.t.}~
\einf_{\bx\in\Omega}q(\bx,p_0)\geq 0.
\end{equation}
Since 
\[
\frac{q(\bx,p)}{p}=-\1/2 \div \ba +\frac{p}{4}|\ba|^2 ,
\]
the inequality
\begin{equation}\label{qp}
\frac{q(\bx,p_1)}{p_1}\geq \frac{q(\bx,p_2)}{p_2}
~~~~~
(p_1\geq p_2>0),
\end{equation}
holds.
Hence, the condition (\ref{condition1}) implies
$\einf_{\bx\in\Omega} q(\bx,p)\geq 0$ for all $p\geq p_0$.
From Lemma~\ref{lemp12}, we have the following theorem.
\begin{Th}
Suppose the condition {\rm (\ref{condition1})}.
Then, $\lambda_1(p)$ is nondecreasing 
for $p\geq p_0/2$.  
In particular, 
the exponential decay phenomenon does not occur
in this case.
\end{Th}
\noindent
{\it Proof.}
We assume that $p_1>p_2\geq p_0/2$.
Then we have $p_1+p_2> p_0$.
From the first inequality of 
Lemma~\ref{lemp12}
and (\ref{qp}),
we obtain
\[
\lambda_1(p_1)-\lambda_1(p_2)
\geq
\frac{p_1-p_2}{p_1+p_2}
\inf q(p_1+p_2)
=
(p_1-p_2)
\inf \frac{q(p_1+p_2)}{p_1+p_2}
\geq
(p_1-p_2)
\inf \frac{q(p_0)}{p_0}
\geq 0.
\]
\qed \\

We notice the following necessary condition for (\ref{condition1}).
\begin{Prop}
The condition
{\rm (\ref{condition1})} implies
that
\begin{equation}\label{condition2}
\min_{\bx\in\ov{\Omega'}}b(\bx)
= \min_{\bx\in \partial \Omega'}b(\bx)
~~~~~(\Omega'
\mbox{\rm : an arbitrary subdomain of $\Omega)$.}
\end{equation}
\end{Prop}
\noindent
{\it Proof.}
We remark that the condition (\ref{condition1})
is equivalent to
$\K b\geq 0$, where $\K$ is a linear differential
operator defined by
\[
\K u:=-\Delta u +\frac{p_0}{2}\nabla b\cdot \nabla u .
\]
Since $b$ is a supersolution with respect to $\K$,
from a consequence of the weak maximum principle
(\cite{G-T83}, Theorem~3.1),
the condition (\ref{condition2})
follows.
\qed \\

An inverse condition of (\ref{condition2})
is given by the following potential well condition.
\begin{Def}\label{pwell}
For a fixed velocity potential $b\in W^{2,\infty}(\Omega)$,
a Lipschitz (nonempty) subdomain $\Omega'\subset \Omega$
is called a potential well if
the condition
\[
\min_{\bx\in\ov{\Omega'}}b(\bx)
< \min_{\bx\in \partial \Omega'}b(\bx),
\]
is satisfied. 
Furthermore,
\[
b_0:=\min_{\bx\in \partial \Omega'}b(\bx)
-\min_{\bx\in\ov{\Omega'}}b(\bx)>0.
\]
is called the depth of a potential well $\Omega'$.
\end{Def}
We remark that
if there exists a potential well,
the condition (\ref{condition2}) does not hold,
and neither does (\ref{condition1}).
Moreover, according to the next theorem,
existence of a potential well 
implies the exponential decay phenomenon
of principal eigenvalues.

Let $\Omega'$ be a subdomain of $C^2$-class
and let $\bnu$ denote the outward unit normal vector
on $\partial \Omega'$.
If it satisfies the condition:
\begin{equation}\label{Fcondition2}
\ba(\bx)\cdot\bnu(\bx)>0
~~~~~(\bx\in\partial\Omega'),
\end{equation}
which is a sufficient
condition for the exponential decay phenomenon in $\Omega'$
obtained in \cite{Fri73},
then the exponential decay phenomenon occurs also in $\Omega$
due to Corollary~\ref{fcor}.
Actually, under the condition (\ref{assumeb}), the condition
(\ref{Fcondition2}) implies that $\Omega'$ is a potential well.

\begin{Th}\label{mth}
We suppose that there exists a potential well $\Omega'$
with depth $b_0>0$.
Then an exponential decay phenomenon occurs:
\[
\liminf_{p\to\infty}\,\frac{1}{p}\log \frac{1}{\lambda_1(p)}
\geq b_0.
\]
In other words,
for any $\omega \in (0,b_0)$,
there exists $C>0$ such that
\[
0<\lambda_1(p)\leq Ce^{-\omega p}
~~~~~
(p\geq 0).
\]
\end{Th}
\noindent
{\it Proof.}
Without loss of generality, we can assume 
$\min_{\sbx\in\ov{\Omega'}}b(\bx) =0$.
From Theorem~\ref{ONS2}, $\lambda_1(p)$ is
given by
\[
\lambda_1(p)=\min_{u\in H^1_0(\Omega),u\not\equiv 0}I_{pb}(u).
\]
Hence, for arbitrary $u\in H^1_0(\Omega)$
$(u\not\equiv 0)$ and $\beta\in\R$,
we have
\begin{equation}\label{luI}
\lambda_1(p)\leq I_{pb}(u)
=
\frac{\int_\Omega e^{-p(b-\beta)}
|\nabla u|^2d\bx}
{\int_\Omega e^{-p(b-\beta)}u^2d\bx}.
\end{equation}

Let us choose $\beta$ satisfying
$0<\beta <\beta+\omega <b_0$.
Since $b\in C^0(\ov{\Omega})$
and $\min_{\partial \Omega'}b=b_0>\beta+\omega$, there exists 
$\vep >0$ such that
$b(\bx)\geq \beta+\omega$ for all
$\bx\in N^{\vep}(\partial \Omega')$,
where
\[
N^{\vep}(\partial \Omega')
:=
\left\{
\bx\in\Omega';~\mbox{\rm dist}(\bx,\partial \Omega')<\vep
\right\},
~~~
\mbox{\rm dist}(\bx,\partial \Omega')
:=\min \{|\bx-\by|;~\by\in \partial \Omega'\}.
\]

We define $\hat{u}\in H^1_0(\Omega)\cap C^{0,1}(\ov{\Omega})$ by
\begin{equation}\label{hatu}
\hat{u} (\bx):=
\left\{
\begin{array}{ll}
0&(\bx\in \ov{\Omega}\setminus \Omega')\\~\\
\vep^{-1}\dist(\bx,\partial \Omega ')&
(\bx\in N^\vep(\partial \Omega'))\\~\\
1&(\bx\in \Omega'\setminus N^\vep(\partial \Omega'))\,,
\end{array}\right. 
\end{equation}
and substitute it to (\ref{luI}).
Since
$|\nabla \hat{u}|=\vep^{-1}$ a.e. in $N^\vep(\partial \Omega')$
and $|\nabla \hat{u}|=0$ in $\Omega\setminus\ov{N^\vep(\partial \Omega')}$,
we obtain
\[
\lambda_1(p)\leq
\frac{\int_{N^\vep(\partial \Omega')}
e^{-p(b(\sbx)-\beta)}
\vep^{-2}d\bx}
{\int_{\Omega'\setminus N^{\vep}(\partial \Omega')}
e^{-p(b(\sbx)-\beta)}d\bx}
\leq
\frac{|N^\vep(\partial \Omega') |\,\vep^{-2}\,
e^{-\omega p}}
{|\{\bx\in\Omega';~b(\bx)\leq \beta\}|}
=
Ce^{-\omega p},
\]
where
\[
C:=\vep^{-2}\,
\left|\{\bx\in\Omega';~b(\bx)\leq \beta\}\right|^{-1}\,
 \left| N^\vep(\partial \Omega') \right|\,.
\]
\qed \\

The above theorem can be extended to 
first $m$ eigenvalues $\lambda_1(p),\cdots,\lambda_m(p)$,
if there are $m$ potential wells.
\begin{Th}\label{mth2}
Let $\Omega_j\subset\Omega$ $(j=1,\cdots,m)$
be disjoint potential wells 
$(\overline{\Omega_i}\cap\overline{\Omega_j}=\emptyset$
for $i\neq j)$, and let $\omega_0$ be
the minimum depth of the potential wells
$\Omega_j$ $(j=1,\cdots,m)$.
Then, for any
$\omega \in (0,\omega_0)$,
there exists $C>0$ such that
\[
0<\lambda_1(p)<\lambda_2(p)\leq \cdots
\leq \lambda_m(p)\leq 
Ce^{-\omega p}
~~~~~
(p\geq 0).
\]
\end{Th}
\noindent
{\it Proof.}
For each $j=1,\cdots,m$, we define 
$\hat{u}_j\in H^1_0(\Omega)\cap C^{0,1}(\ov{\Omega})$
with $\supp (\hat{u}_j)=\ov{\Omega_j}$
in the same manner as (\ref{hatu}).
Let X be a linear subspace of $H^1_0(\Omega)$
generated by functions
$\widehat{u}_1,\cdots,\widehat{u}_m$. 
It is easy to see $\mathrm{dim}\, X=m$.
Due to the min-max formula in Theorem~\ref{ONS2},
the $m$th eigenvalue $\lambda_m(p)$ is estimated from above as
\begin{equation}\label{minmax}
\lambda_m(p)\leq 
\max_{u\in X\setminus \{0\}} I_{pb}(u)
=
\max_{c_j} \,
I_{pb}(c_1 \hat{u}_1+\cdots + c_m \hat{u}_m),
\end{equation}
where in $\max_{c_j}$, the maximum is taken over
all $(c_1,\cdots,c_m)\in\R^m\setminus \{(0,\cdots,0)\}$.
Since $\ov{\Omega_j}$ are disjoint, we have
\begin{eqnarray*}
I_{pb}(c_1 \hat{u}_1+\cdots + c_m \hat{u}_m)
&=&
\frac{\DS{
\int_\Omega e^{-pb(\sbx)}
\left(\sum_{j=1}^m c_j^2 |\nabla \hat{u}_j (\bx)|^2\right) d\bx}}
{\DS{\int_\Omega e^{-pb(\sbx)}
\left(\sum_{j=1}^m c_j^2 |\hat{u}_j (\bx)|^2\right) d\bx}}=
\frac{\DS{\sum_{j=1}^m c_j^2
\int_{\Omega_j} e^{-pb(\sbx)}|\nabla \hat{u}_j (\bx)|^2 d\bx}}
{\DS \sum_{j=1}^m c_j^2 \int_{\Omega_j} e^{-pb(\sbx)}|\hat{u}_j (\bx)|^2 d\bx}\\~\\
&=&
\frac{\DS \sum_{j=1}^m d_j^2\ 
\frac{\DS \int_{\Omega_j} e^{-pb(\sbx)}|\nabla \hat{u}_j (\bx)|^2 d\bx}
{\DS \int_{\Omega_j} e^{-pb(\sbx)}|\hat{u}_j (\bx)|^2 d\bx}}
{\DS \sum_{j=1}^m d_j^2}\leq 
\max_{1\leq j\leq m}
\frac{\DS \int_{\Omega_j} e^{-pb(\sbx)}|\nabla \hat{u}_j (\bx)|^2 d\bx}
{\DS \int_{\Omega_j} e^{-pb(\sbx)}|\hat{u}_j (\bx)|^2 d\bx},
\end{eqnarray*}
where
we have defined $d_j:=c_j 
(\int_{\Omega_j} e^{-pb(\sbx)}|\hat{u}_j (\bx)|^2 d\bx)^{1/2}$.

Similar to the proof of Theorem~\ref{mth},
choosing sufficiently small $\vep>0$, we can show that
the last term is bounded by $Ce^{-\omega p}$ from above.
\qed 

\subsection{\large Precise asymptotic behaviour in one dimensional case}
\label{subsec2} 
\setcounter{equation}{0}
In this section,
we introduce a more precise estimate for the 
exponential decay phenomenon of 
the principal eigenvalues in one dimensional case.
In \cite{K-N05},  some of results in this section 
have been already announced in Japanese by one of the authors.

We consider the following one dimensional eigenvalue problem
on $\Omega=(-l,l)$, where $l>0$.
\begin{equation}\label{1d-problem}
\left\{
\begin{array}{ll}
-u''(x)+pa(x)u'(x)=\lambda u (x)
&(-l<x<l)\\~\\
u (-l)=u (l)=0\\~\\
u (x)>0&(-l<x<l).
\end{array}
\right.
\end{equation}
The principal eigenvalue is denoted by $\lambda_1(p)>0$.
In the following argument,
we fix one of the primitive functions
of $a\in L^\infty(-l,l)$ and denote it by 
$b\in W^{1,\infty}(-l,l)=C^{0,1}([-l,l])$, 
i.e., 
\[
b'(x)=a(x)~~~\mbox{a.e.}~~x\in (-l,l).
\]
We also define
\[
b_0:=b_2-b_1,~~~~~
b_1:= \min_{0 \leq x \leq l} b(x),~~~~~
b_2:=\max_{0 \leq x \leq l} b(x),
\]
\[
B_i:=\{x\in [0,l];~b(x)=b_i\}
~~~~~(i=1,\,2).
\]
We assume the following
assumption.
\begin{Assum}\label{AB}
The function $a(x)$ is an odd function 
belonging to $L^\infty (-l,l)$, and it satisfies
\[
\max_{x\in B_1} x <\min_{y\in B_2} y\,.
\]
\end{Assum}
\noindent
A typical example of $a(x)$ 
satisfying Assumption~\ref{AB} is
\begin{equation}\label{a-ex1}
a(x)=\frac{|x|^\alpha}{x}~~~~~
b(x)=\frac{1}{\alpha}|x|^\alpha,
\end{equation}
for $\alpha \geq 1$ and $l>0$.
In case that
\begin{equation}\label{a-ex2}
a (x)=\sin x,~~~~~
b(x)=-\cos x,
\end{equation}
it satisfies Assumption~\ref{AB} if $0<l<2\pi$.
Another example is
\begin{equation}\label{a-ex3}
a (x)=x^3-x,~~~~~
b(x)=\frac{1}{4}(x^2-1)^2,
\end{equation}
which satisfies Assumption~\ref{AB} if $l>\sqrt{2}$.

In the following arguments, $A(p)\sim B(p)$ as $p\to\infty$
means that $A(p)/B(p)\to 1$ as $p\to\infty$.
Now we present the main result of this section.
\begin{Th}\label{1dth}
Under Assumption~\ref{AB}, we have
\begin{equation}\label{assertion}
\lambda_1(p)\sim 
\left( \int_0^l e^{-pb(x)} dx \right)^{-1} 
\left( \int_0^l e^{pb(y)} dy \right)^{-1}
~~~~~\mbox{as}~~p \to \infty .
\end{equation}
\end{Th}
\noindent
{\it Proof.}
We first remark that
the principal eigenfunction $u_1(x)$ defined by
(\ref{1d-problem}) is an even function since
$a(x)$ is odd. 
By means of the change of variables from $x\in [0,l]$
to $r\in [0,1]$ by $r=\rho(x)$:
\[
\rho (x):=\frac{1}{\sigma} \int_0^x e^{pb(y)}dy~~~(0\le x \le l),
~~~~~
\sigma :=\int_0^l e^{pb(x)}dx,
\]
\[
U_1(r):=u_1(\rho^{-1}(r))~~~~~(0\leq r\leq 1),
\]
it follows that $U_1(r)$ and $\lambda_1(p)$ satisfy
\begin{equation}\label{Uprob}
\left\{
\begin{array}{ll}
-\,U_1''(r)=\lambda_1(p)\,z(r)\,U_1(r) &(0<r<1)\\~\\
U_1'(0)=0,~~~U_1(1)=0&  \\~\\
U_1(r)>0 & (0<r<1),
\end{array}
\right.
\end{equation}
where 
\[
z(r):= \sigma ^2 e^{-2p\,b \left( \rho ^{-1} (r) \right)}>0
~~~~~(0\le r \le 1).
\]
We remark that the quantity in the right hand side of
(\ref{assertion}) is given by 
the following equality
\[
\int_0^1 z(r)dr 
=
\left( \int_0^l e^{-pb(x)} dx \right)
\left( \int_0^l e^{pb(y)} dy \right) ,
\]
which can be checked by direct calculation.

The eigenvalue problem (\ref{Uprob}) is selfadjoint and 
it is found that $\lambda_1(p)$ is given by 
minimum of the Rayleigh quotient:
\begin{equation}\label{l1p}
\lambda_1(p)
=
\frac{\int_0^1 |U_1'(r)|^2 dr}
{\int_0^1 z(r) U_1 (r)^2dr}
= \min _{\varphi \in V} 
\frac{\int_0^1 |\varphi '(r)|^2dr}{\int_0^1 z(r) \varphi (r)^2dr},
\end{equation}
where
$V:= \{ \varphi \in H^1(0,1);\,
\varphi (1)=0,~~\varphi \not\equiv 0 \}$.

Since
\[
U_1'(r)=U_1'(0)+\int_0^r U_1''(s)ds =
-\lambda_1(p) \int_0^r z(s)U_1(s) ds <0
~~~~~
(0<r\leq 1),
\] 
we have $\max U_1=U_1(0)$.
Without loss of generality,
we can assume that $\max U_1=U_1(0)=1$.
Among functions $\varphi\in H^1(0,1)$ with 
$\varphi(0)=1$ and $\varphi(1)=0$, the minimum 
of the Dirichlet integral $\int_0^1 |\varphi'(r)|^2 dr$ 
is achieved by $\varphi(r)=1-r$.
Hence, from the first equality of (\ref{l1p}),
we get an estimate from below:
\begin{equation}\label{ebelow}
\lambda_1 (p) = 
\frac{\int_0^1 | U_1 '(r)|^2 dr}{\int_0^1 z(r) U_1 (r)^2dr}
\geq
\frac{\int_0^1 |(1-r)'|^2 dr}{\int_0^1 z(r) dr}
=
\frac{1}{\int_0^1 z(r)dr}\,.
\end{equation}

On the other hand,
from the second equality of (\ref{l1p}),
substituting $\varphi(r)=1-r$,
we obtain an estimate from above:
\begin{equation}\label{eabove}
\lambda_1 (p) 
=
\min_{\varphi \in V} 
\frac{\int_0^1 | \varphi '(r)|^2 dr}{\int_0^1 z(r) \varphi (r)^2dr}
\leq
\frac{\int_0^1 | (1-r)'|^2 dr}{\int_0^1 z(r) (1-r)^2 dr}
= 
\frac{1}{\int_0^1 z(r) (1-r)^2 dr}.
\end{equation}
From these two estimates (\ref{ebelow}) and (\ref{eabove}), 
if we can show that
\begin{equation}\label{limit}
\frac{\int_0^1 z(r) (1-r)^2 dr}{\int_0^1 z(r)dr}
\to 1
~~~~~\mbox{as}~~p\to\infty,
\end{equation}
then we obtain the assertion of the theorem.

Let us show (\ref{limit}).
We define
\[
Z_m(p):=\int_0^1 z(r)r^mdr >0
~~~~~(m=0,\,1,\,2).
\]
Then we have
\[
\frac{\int_0^1 z(r) (1-r)^2 dr}{\int_0^1 z(r)dr}
=
1-\frac{2Z_1(p)-Z_2(p)}{Z_0(p)}.
\]
Since $0\leq Z_2(p)\leq Z_1(p)$ holds,
it is sufficient to show that $Z_1(p)/Z_0(p)\to 0$ as
$p\to \infty$.

For $m=0,~1$, we have the following equalities:
\[
Z_m(p)=\sigma \int_0^l e^{-pb(x)}\rho (x)^mdx
=\int\mbox{\hspace{-1ex}}\int_{D_m}e^{p(b(y)-b(x))}dxdy,
\]
where
\[
D_0:= \{ (x,y)~|~0<x<l,~0<y<l \},
~~~~~
D_1:= \{ (x,y)~|~0<y<x<1\} .
\]
We define
\[
\omega_m:=\max_{(x,y)\in \overline{D_m}} (b(y)-b(x))~~~~~(m=0,1).
\]
Then, $\omega_0=b_0$ and Assumption~\ref{AB} is equivalent to 
$\omega_0>\omega_1$.
Defining 
\[
G:=\left\{(x,y)\in D_0~|~b(y)-b(x)\geq \frac{\omega_0+\omega_1}{2}\right\},
\]
we obtain $|G|>0$ and 
\begin{equation}\label{I_2}
\frac{Z_1(p)}{Z_0(p)}=
\frac{\int\mbox{\hspace{-1ex}}\int_{D_1}e^{p(b(y)-b(x)-\omega_1)}dxdy}
{\int\mbox{\hspace{-1ex}}\int_{D_0}e^{p(b(y)-b(x)-\omega_1)}dxdy}
\leq
\frac{\int\mbox{\hspace{-1ex}}\int_{D_1}dxdy}
{\int\mbox{\hspace{-1ex}}\int_{G}e^{p(\frac{\omega_0+\omega_1}{2}-\omega_1)}dxdy} 
=\frac{|D_1|}{|G|}e^{-p\frac{\omega_0-\omega_1}{2}}\to 0, 
\end{equation}
as $p$ tends to infinity.
Hence we have proved (\ref{limit}).
\qed \\

As we have seen in the above proof, 
the exponential decay phenomenon of 
the principal eigenvalue occurs
under the Assumption~\ref{AB}.
\begin{Cor}
Under Assumption~\ref{AB}, we have
\begin{equation} \label{eq:Friedman}
\lim_{p \to \infty} \,\frac{1}{p} \log \frac{1}{\lambda_1 (p)}  =b_0.
\end{equation}
In other words, for arbitrary $\vep\in (0,b_0)$,
there exists $p_0>0$ such that
\begin{equation} \label{cor_lambda}
e^{-p(b_0+ \vep)} \leq \lambda_1 (p) 
\leq e^{-p(b_0- \vep)}~~~~~(p \geq p_0).
\end{equation}
\end{Cor}
\noindent
{\it Proof.}
It is sufficient to show
(\ref{eq:Friedman}), since (\ref{cor_lambda}) is
equivalent to (\ref{eq:Friedman}).
From Theorem~\ref{1dth}, 
$\lambda_1(p) Z_0(p)\to 1$ as $p\to\infty$
holds.
We have (\ref{eq:Friedman}) by the equalities:
\begin{equation}\label{lp}
\lim_{p\to\infty}
\left( \frac{1}{p} \log \frac{1}{\lambda_1 (p)} -b_0\right)
=
\lim_{p\to\infty}
\left(\frac{1}{p} \log Z_0 (p) -b_0\right)
=
\lim_{p\to\infty}
\frac{1}{p} \log \left(e^{-pb_0}Z_0 (p)\right)
=0,
\end{equation}
where the last equality is shown as follows.
From the inequality $e^{-pb_0}Z_0 (p)\leq |D_0|$, 
it follows that, 
\[
\limsup_{p\to\infty}
\frac{1}{p} \log \left(e^{-pb_0}Z_0 (p)\right)
\leq
\lim_{p\to\infty}
\frac{\log |D_0|}{p}
=0.
\]
Moreover, for arbitrary $\vep >0$, we define
$D_0^\vep :=\left\{
(x,y)\in D_0;~b(y)-b(x)\geq b_0-\vep \right\}$.
Then, from $|D_0^\vep |>0$ and the inequality:
\[
e^{-pb_0}Z_0 (p)
=
\int\mbox{\hspace{-1ex}}\int_{D_0}
e^{p(b(y)-b(x)-b_0)}dxdy 
\geq 
\int\mbox{\hspace{-1ex}}\int_{D_0^\vep}
e^{-p\vep}dxdy
=
e^{-p\vep}|D_0^\vep | ,
\]
we obtain
\[
\liminf_{p\to\infty}
\frac{1}{p} \log \left(e^{-pb_0}Z_0 (p)\right)
\geq
\lim_{p\to\infty}
\frac{-p\vep +\log |D_0^\vep |}{p}
=
-\vep .
\]
Since $\vep >0$ is arbitrary, the last equality of (\ref{lp})
follows.
\qed \\

For concrete examples of $a(x)$,
more precise 
asymptotic behaviours of $\lambda_1(p)$ 
can be derived from Theorem~\ref{1dth}.
For example, if $|B_1| |B_2|>0$ then we get
\[
\lambda_1(p)
\sim \frac{1}{|B_1| |B_2|}
e^{-pb_0}
~~~\mbox{as}~p\to\infty ,
\]
since the following estimates hold:
\[
\int_0^le^{-pb(x)}dx\sim |B_1|e^{-pb_1}
~~~\mbox{as}~p\to\infty ,~~~
\mbox{if}~|B_1|>0, 
\]
\[
\int_0^le^{pb(x)}dx\sim |B_2|e^{pb_2}
~~~\mbox{as}~p\to\infty ,
~~~\mbox{if}~|B_2|>0.
\]

In case that $|B_1| |B_2|=0$, the following lemma is useful.
\begin{Lem}\label{lem-alpha}
Let $L>0$ and $\mu >0$.
Suppose that $g\in C^0([0,L])$ satisfies $g(x)>0$ for $x\in (0,L]$
and $\lim_{x\downarrow 0}x^{-\mu}g(x)=1$.
Then we have
\[
\int_0^L e^{-pg(x) }dx
\sim \Gamma \left( \frac{1}{\mu}+1\right) p^{-\frac{1}{\mu}}
~~~\mbox{as}~p\to\infty ,
\]
where $\Gamma$ stands for the Gamma function 
$\Gamma (\zeta)=\int_0^\infty s^{\zeta-1}e^{-s}ds$.
\end{Lem}
\noindent
{\it Proof.}
In case that $g(x)=x^\mu$, by the change of variables $s=px^\mu$,
we have
\[
\frac{p^{\frac{1}{\mu}}}{\Gamma (1/\mu +1 )} 
\int_0^L e^{-px^\mu }dx
=
\frac{p^{\frac{1}{\mu}}}{\Gamma (1/\mu +1 )} 
\frac{1}{\mu p^{\frac{1}{\mu}}}
\int_0^{pL^\mu} s^{\frac{1}{\mu}-1}e^{-s}ds
=
\frac{1}{\Gamma (1/\mu)} 
\int_0^{pL^\mu} s^{\frac{1}{\mu}-1}e^{-s}ds
\to
1,
\]
as $p$ tends to infinity. Hence, the assertion follows 
for $g(x)=x^\mu$. For general $g(x)$, we omit a proof but
the same asymptotic behaviour can be shown.
\qed \\

By using this lemma, for example, we obtain the following proposition.
\begin{Prop}
For the case {\rm (\ref{a-ex1})} with $\alpha \geq 1$, we have
\[ 
\lambda_1 (p) \sim 
\frac{\left(\frac{1}{\alpha}\right)^{\frac{1}{\alpha}}l^{\alpha -1}}
{\Gamma \left( \frac{1}{\alpha}+1\right)}
\,p^{\frac{1}{\alpha}+1}\,
e^{-\frac{l^\alpha}{\alpha}p}
~~~\mbox{as}~p\to\infty .
\]
In particular, if $a(x)=x$ 
$(\,\mbox{i.e.},
-u''(x)+pxu'(x)= \lambda u (x)~~(-l<x<l))$,
then
\begin{equation}\label{a=x}
\lambda_1 (p) \sim \sqrt{\frac{2}{\pi}}\, l\,
p^{\frac{3}{2}}\, e^{-\frac{l^2}{2}p}
~~~\mbox{as}~p\to\infty .
\end{equation}
\end{Prop}
\noindent
{\it Proof.}
From Lemma~\ref{lem-alpha}, we have
\[
\int_0^l e^{-pb(x) }dx
=
\int_0^l e^{-\frac{p}{\alpha}x^\alpha }dx
\sim \Gamma \left( \frac{1}{\alpha}+1\right) 
\left(\frac{p}{\alpha}\right)^{-\frac{1}{\alpha}}
~~~\mbox{as}~p\to\infty .
\]
Using the change of variables $y=l-x$, we also have 
\[
\int_0^l e^{pb(x) }dx
=
\int_0^l e^{\frac{p}{\alpha}x^\alpha }dx
=
e^{\frac{p}{\alpha}l^\alpha}
\int_0^l e^{-l^{\alpha -1}pg(y)} dy
\sim 
e^{\frac{p}{\alpha}l^\alpha}
\Gamma (2)
\left(l^{\alpha -1}p\right)^{-1}
~~~\mbox{as}~p\to\infty ,
\]
where $g(y):=\alpha^{-1}l^{1-\alpha}(l^\alpha -(l-y)^\alpha)$
satisfies the condition of Lemma~\ref{lem-alpha} with $\mu=1$.
From Theorem~\ref{1dth}, we obtain the assertion.
\qed \\

Similarly, we can calculate that, for (\ref{a-ex2}):
\[
\lambda_1(p)\sim
\left\{
\begin{array}{ll}
\frac{\sqrt{2}\sin l}{\sqrt{\pi}}
\,p^{\frac{3}{2}}\,e^{-(1-\cos l)p}
&(\mbox{if}~0<l<\pi)\\~\\
\DS{\frac{1}{2\pi}\,p\, e^{-2p}}
&(\mbox{if}~l=\pi)\\~\\
\DS{\frac{1}{\pi}\,p\, e^{-2p}}
&(\mbox{if}~\pi<l<2\pi)
\end{array}\right.
~~~\mbox{as}~p\to\infty ,
\]
and 
for (\ref{a-ex3}):
\[
\lambda_1(p)\sim
\frac{a(l)}{\sqrt{\pi}}\,p^{\frac{3}{2}}\, e^{-b(l)p}
~~~\mbox{as}~p\to\infty ,
~~\mbox{if}~l>\sqrt{2} .
\]

We remark that similar precise asymptotic expansion 
has been studied in de Groen \cite{deG80a} by means of 
precise approximation of the eigenfunctions
under a different assumption.
Our result (\ref{a=x}) coincides to the result obtained 
there ((8.11) in \cite{deG80a}).

At the end of this section,
we remark that Theorem~\ref{1dth} can not be 
straightforwardly extended
to multi-dimensional cases.
To see this, we assume that $\Omega$ and $a(\bx)$ satisfy the conditions of
Proposition~\ref{rect} with $\cI_i=(-l_i,l_i)$ and
that each $a_i(x)$ satisfies Assumption~\ref{AB}.
Then, since $b(\bx)=\sum_{i=1}^nb_i(x_i)$ with $b_i'(x_i)=a_i(x_i)$,
we have
\[
\left( \int_\Omega e^{-pb(\sbx)} d\bx \right)^{-1} 
\left( \int_\Omega e^{pb(\sby)} d\by \right)^{-1}
\sim
4^{-n}\prod_{i=1}^n \lambda_1^{(i)}(p)
~~~\mbox{as}~p\to\infty ,
\]
whereas, from Proposition~\ref{rect}, we obtain
\[
\lambda_1(p)=\sum_{i=1}^n \lambda_1^{(i)}(p).
\]
In this case, the asymptotic behaviour of $\lambda_1(p)$
is determined by the slowest decaying $\lambda_1^{(i)}(p)$.

\subsection{\large Chemotaxis effect in biological colonies}
\label{chemotaxis} 
\setcounter{equation}{0}

A simple biological interpretation
of the exponential decay phenomenon of principal
eigenvalue is discussed in this section.
We consider the following simple PDE model
for aggregation phenomena in biology.

Let $v (\bx,t)\geq 0$ be a population density at $\bx\in\Omega$
of an organism (e.g.\,bacteria), 
which exhibits an aggregation behaviour
due to an attractive chemical substance.
This property is called chemotaxis and 
the attractive chemical substance is called chemotaxis substance.
We denote by $c(\bx,t)\geq 0$
the density of the chemotaxis substance.
Then one of the standard PDE description of the chemotaxis 
effect is given by the following equation.
\begin{equation}\label{cxm}
v_t=\Delta v - p \,\div (v \nabla c)
~~~~~\mbox{in}~\Omega\times (0,\infty),
\end{equation}
where $ p >0$ stands for 
a chemotaxis parameter. If $ p $ is larger,
$v $ exhibits stronger chemotaxis.

Keller and Segel (\cite{K-S70}, \cite{K-S71})
and other great number of literatures
(for example Diaz and Nagai \cite{D-N95}, Stevens \cite{St00}
etc.)
have proposed (\ref{cxm}) nonlinearly coupled with
the production of chemotaxis substance by $v $:
\[
\alpha c_t= \Delta c -\gamma c + v 
~~~~~\mbox{in}~\Omega\times (0,\infty),
\]
with a suitable boundary condition for $c$,
where $\alpha$ and $\gamma$ are nonnegative constants.
There are several results on finite time blow-up
of these systems if the initial density $v (\bx,0)$ is
large enough. 
The blowup of the solution in the Keller-Segel
model has attracted many mathematicians' interests
and its mathematical mechanism has been studied.
Please see the above literature
and references therein for more details.
Especially, a good short review for mathematical results
on the Keller-Segel model is found in
\S~8 of \cite{St00}.

On the other hand, in \cite{D-N95}, 
under the zero Dirichlet boundary condition for $v $
with $\alpha=0$,
they proved that $v $ decays exponentially to zero 
as $t\to\infty$ if $n=1$ or if $v (\bx,0)$ is
sufficiently small for $n\geq 2$.
The zero Dirichlet boundary condition for $v $
biologically means the extinction of the organism
at the boundary. 
The result in \cite{D-N95}
can be interpreted as that,
the colony formed by the chemotaxis aggregation 
can not sustain eternally against the diffusion and
the boundary extinction effect 
($v=0$ on $\partial \Omega$), if the colony size is not
large enough.

Let us consider a simpler linear version of
this chemotaxis aggregation 
vs boundary extinction problem.
We suppose that the chemotaxis substance $c$
is a priori given and does not change in time, i.e. $c=c(\bx)$,
and we consider 
the initial-boundary value problem:
\begin{equation}\label{linearchm}
\left\{
\begin{array}{ll}
v_t=\Delta v - p \,\div (v \nabla c(\bx))
&\mbox{in}~\Omega\times (0,\infty)\\
~&\\
v =0
&\mbox{on}~\partial\Omega\times (0,\infty)\\
~&\\
v (\cdot,0)=v_0(\bx)
&\mbox{in}~\Omega ,
\end{array}
\right.
\end{equation}
where $v_0(\bx)\geq 0$ is a given initial
density function with $v_0\not\equiv 0$.
It is mathematically well-known that
the boundary extinction is always 
stronger than the chemotaxis aggregation.
It decays exponentially as
\[
v (\bx,t)\sim C e^{-\lambda_1( p ) \,t}v_1(\bx,p)
~~~~~
(t\to\infty),
\]
where $\lambda_1( p )$ is nothing but
$\lambda_1(p)$ in Section~\ref{asymptotic} with 
$b(\bx)=-c(\bx)$,
and $v_1(\bx,p)$ is the principal eigenfunction of $K_{pb}^*$.
If $b(\bx)=-c(\bx)$ has the potential well (Definition~\ref{pwell})
of depth $\omega >0$, $\lambda_1( p )$ exhibits 
the exponential decay phenomenon:
$\lambda_1( p )\approx O(e^{-\omega  p })$ for large $ p $.
The potential well condition for $-c(\bx)$
means that there is a concentration peak of 
the chemotaxis substance of
hight $\omega$.

The reciprocal of the principal eigenvalue 
$1/\lambda_1( p )$ stands for 
the life span of the colony,
since the half-life of exponential decaying quantity
$e^{-\lambda t}$
is given by $(\log 2)/\lambda$.
Theorem~\ref{mth} is interpreted as that
if $c(\bx)$ has a concentration peak
of hight $\omega >0$ then
the life span of the colony is
exponentially long as
\[
\frac{1}{\lambda_1( p )}\approx 
Ce^{\omega  p }
~~~~~
\mbox{as}~~ p \to\infty .
\]
In other words, the colony can not sustain
eternally, but
it sustains in enough long time $O(e^{\omega  p })$
if the chemotaxis parameter is sufficiently large.

\begin{figure}[tbh]
\begin{center}
\begin{minipage}{0.45\hsize}
\begin{center}
\rotatebox{270}{\includegraphics[width=0.75\hsize,clip]{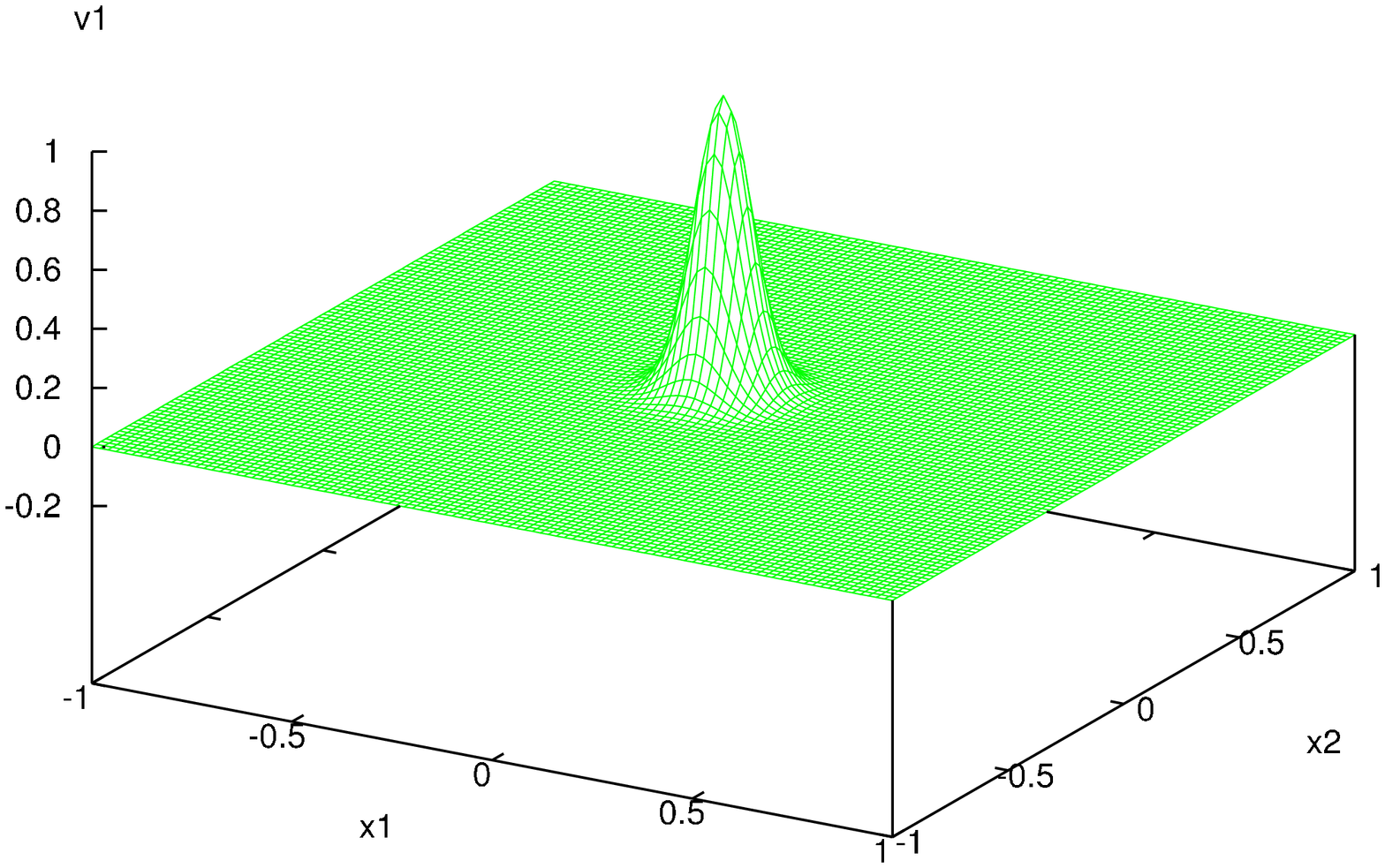}}
\caption{Eigenfunction $v_1$ of (\ref{eigen-u}) 
corresponding to $u_1$ of Figure~\ref{Sim1-profile-p40}.}
\label{Sim1_t=1s_p=40_-pb_3d_c}
\end{center}
\end{minipage}\\
\begin{minipage}{0.45\hsize}
\begin{center}
\rotatebox{270}{\includegraphics[width=0.75\hsize,clip]
{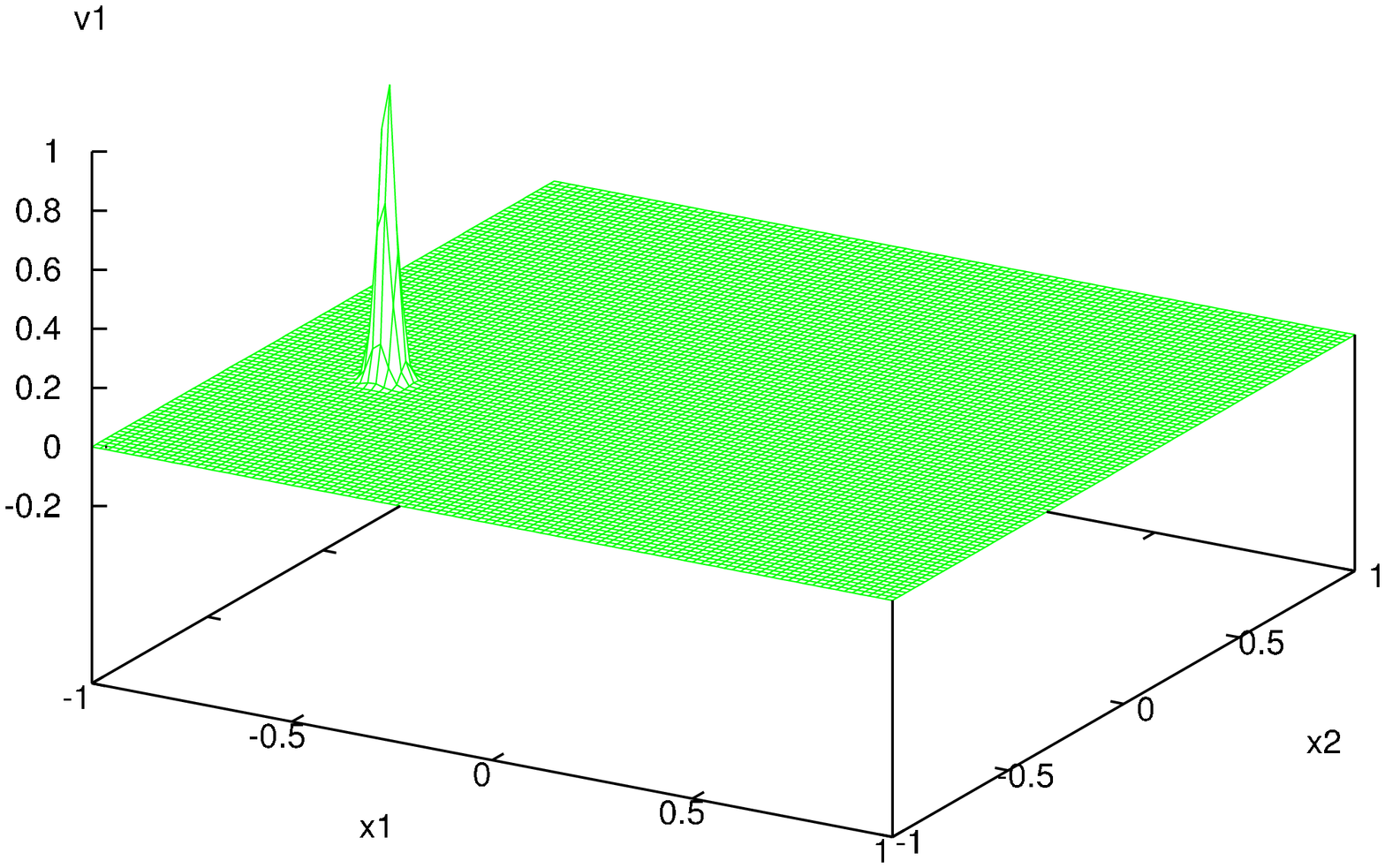}}
\end{center}
\end{minipage}
\hspace{0.08\hsize}
\begin{minipage}{0.45\hsize}
\begin{center}
\rotatebox{270}{\includegraphics[width=0.75\hsize,clip]
{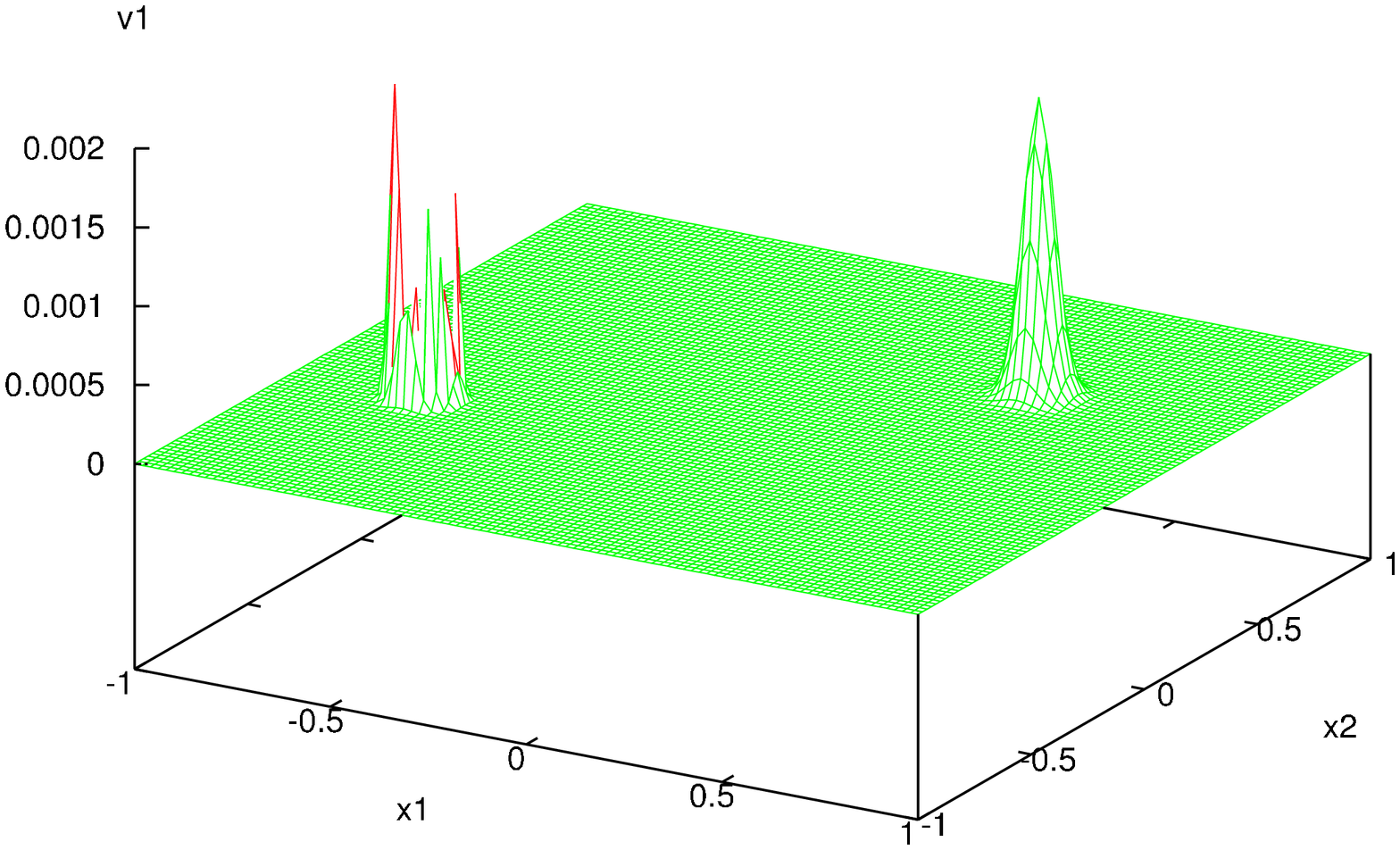}}
\end{center}
\end{minipage}
\caption{Eigenfunction $v_1$ of (\ref{eigen-u}) 
corresponding to $u_1$ of Figure~\ref{Sim3-profile-p100} 
(left), and a scaled figure enlarged in the vertical axis (right).}
\label{Sim3_t=1s_p=100_-pb_3d_c}
\end{center}
\end{figure}

In Figure~\ref{Sim1_t=1s_p=40_-pb_3d_c},
a typical profile of colony is shown. 
This is the eigenfunction $v_1$ of (\ref{eigen-u}) 
with same $\ba(\bx)\,(=-\nabla c(\bx))$ as Example~\ref{Sim1} and $p=40$.
Although
the hight of this colony is decreasing in time according to the linear law,
but the decay is exponentially slow for large $p$
and it substantially forms a very solid colony.
Another profile of $v$ corresponding to Example~\ref{Sim3} with $p=100$
is shown in Figure~\ref{Sim3_t=1s_p=100_-pb_3d_c},
which has two colonies. The hight of the right small colony is
approximately $1/500$ of the hight of the left one.
This is interesting in comparison with Figure~\ref{Sim3-profile-p100}
where the eigenfunction $u_1$ has two terraces with almost same hight.

~\\
\noindent
{\bf Acknowledgement}
We would like to express our thanks to 
Danielle~Hilhorst (Universit\'{e} Paris-Sud),
Kazuhiro~Kurata (Tokyo Metropolitan University)
and Toshitaka~Nagai (Hiroshima University)
for their fruitful discussions and kind suggestions of literature.
The research of S.~Jimbo was partially supported by 
JSPS Grant-in-Aid for Scientific Research (B), No.17340042.
The research of M.~Kimura was partially supported by the Ne\v{c}as Center 
for Mathematical Modelling LC06052 financed by MSMT
while his stay in Czech Technical University in Prague.
The research of H.~Notsu was partially supported by 
JSPS Grant-in-Aid for Scientific Research (S), No.16104001.

\end{document}